\documentclass[a4paper,10pt]{amsart}
\usepackage{amsmath}
\usepackage{amssymb}
\usepackage{amsfonts}
\usepackage{mathtools}
\usepackage{amscd}
\usepackage{todonotes}
\usepackage{marginnote}
\usepackage{latexsym,graphicx}
\usepackage[latin1]{inputenc}
\usepackage{xspace}
\usepackage{xcolor}
\usepackage{array}
\usepackage{newclude}
\usepackage{DotArrow}
\usepackage{tensor}
\usepackage[all]{xy}
\usepackage{tikz}
\usepackage{tikz-cd}
\usepackage{enumitem,kantlipsum}
\usepackage{newclude}
\usepackage[backref=page, bookmarksopen=true]{hyperref}
\hypersetup{pdfpagelabels,
                  plainpages=false,
                    colorlinks=true,       
                    linkcolor=blue}
\renewcommand*{\HyperDestNameFilter}[1]{\jobname-#1} 
\numberwithin{equation}{section}
\usepackage[a4paper, total={5in, 9in}, centering]{geometry}
\usepackage{todonotes}
\usepackage{mathrsfs}

\usepackage[nameinlink]{cleveref}

\makeatletter
\providecommand{\leftsquigarrow}{%
  \mathrel{\mathpalette\reflect@squig\relax}%
}
\newcommand{\reflect@squig}[2]{%
  \reflectbox{$\m@th#1\rightsquigarrow$}%
}

\newcommand{\sspace}{\vspace{0.25cm}}
\newcommand{\noi}{\noindent}

\theoremstyle{plain}
\newtheorem{theor}{Theorem}[section]
\newtheorem{conj}[theor]{Conjecture}
\newtheorem{prop}[theor]{Proposition}
\newtheorem{lem}[theor]{Lemma}

\theoremstyle{remark}
\newtheorem{rem}[theor]{Remark}

\newtheorem{Example}[theor]{Example}

\newtheorem{Question}[theor]{Question}

\theoremstyle{plain}

\newtheorem{defi}[theor]{Definition}

\newcommand{\CC}{{\mathbb C}}
\newcommand{\RR}{{\mathbb R}}

\newcommand{\QQ}{{\mathbb Q}}

\newcommand{\ZZ}{{\mathbb Z}}

\newcommand{\BB}{{\mathbf B}}
\newcommand{\G}{{\mathbf G}}

\newcommand{\NN}{{\mathbb N}}

\newcommand{\Qbar}{{\overline{\QQ}}}

\newcommand{\Spec}{{\rm Spec}}

\newcommand{\Pic}{{\rm Pic}}

\newcommand{\Gr}{{\textnormal{Gr}}}

\newcommand{\ti}[1]{\mbox{$\tilde{#1} $}}
\newcommand{\wti}[1]{\mbox{$\widetilde{#1} $}}

\newcommand{\ol}{\overline}

\newcommand{\End}{{\rm End}}

\newcommand{\Res}{{\rm Res}}

\newcommand{\Gal}{{\rm Gal}}

\newcommand{\Gm}{{\mathbf{G}_{m}}}

\newcommand{\GL}{{\rm \bf GL}}

\newcommand{\tr}{\textnormal{tr}}

\newcommand{\proj}{{\mathbf P}}

\newcommand{\Sp}{\mathbf{Sp}}
\newcommand{\Aut}{\textnormal{Aut}}
\newcommand{\an}{\textnormal{an}}

\newcommand{\HHom}{\textnormal{Hom}}

\newcommand{\pure}{\textnormal{p}}

\def\fH{\mathfrak{H}}

\def\fA{\mathfrak{A}}

\def\fg{\mathfrak{g}}
\def\fh{\mathfrak{h}}

\def\fb{\mathfrak{b}}

\newcommand{\cM}{{\mathcal M}}

\newcommand{\cT}{{\mathcal T}}
\newcommand{\cS}{{\mathcal S}}
\newcommand{\cC}{{\mathcal C}}
\newcommand{\cE}{{\mathcal E}}

\newcommand{\cA}{{\mathcal A}}

\newcommand{\cO}{{\mathcal{O}}}

\newcommand{\cV}{{\mathcal V}}
\newcommand{\cY}{{\mathcal Y}}
\newcommand{\cJ}{{\mathcal J}}
\newcommand{\cU}{{\mathcal U}}

\newcommand{\nH}{\textnormal{H}}

\newcommand{\nS}{\textnormal{S}}
\newcommand{\nM}{\textnormal{M}}
\newcommand{\nV}{\textnormal{V}}

\newcommand{\bH}{{\mathbb H}}

\newcommand{\bV}{{\mathbb V}}

\newcommand{\bL}{{\mathbb L}}
\newcommand{\bM}{{\mathbb M}}

\newcommand{\bfG}{{\mathbf G}}
\newcommand{\bfH}{{\mathbf H}}
\newcommand{\bfU}{{\mathbf U}}
\newcommand{\bfK}{{\mathbf K}}

\newcommand{\oQ}{\overline{\QQ}}

\newcommand{\codim}{\textnormal{codim}}

\newcommand{\Ext}{\textnormal{Ext}}

\newcommand{\Alb}{{\rm Alb}}
\newcommand{\PGL}{{\mathbf{PGL}}}

\newcommand{\SL}{{\mathbf{SL}}}

\newcommand{\Lie}{{\rm Lie}}

\newcommand{\Zar}{\textnormal{Zar}}

\newcommand{\red}{\textnormal{red}}

\newcommand{\MHS}{\textnormal{MHS}}
\newcommand{\VMHS}{\textnormal{VMHS}}
\newcommand{\Betti}{\textnormal{Betti}}
\newcommand{\alb}{\textnormal{alb}}
\newcommand{\PV}{\proj \textnormal{V}}
\newcommand{\Div}{\textnormal{div}}

\newcommand{\pr}{\textnormal{pr}}
\newcommand{\rel}{\textnormal{rel}}
\newcommand{\Jac}{\textnormal{Jac}\,}

\begin{document}

\title[Abelian differentials and their periods]{Abelian differentials and their periods: the bi-algebraic point
of view}
\author{Bruno Klingler}
\address{Bruno Klingler: Humboldt Universit\"at zu Berlin} 
\email{bruno.klingler@hu-berlin.de}
\author{Leonardo A. Lerer}
\address{Leonardo A. Lerer: Weizmann Institute of Science} 
\email{a.leonardo.lerer@gmail.com}

\begin{abstract}
We study the transcendence of periods
of abelian differentials, both at the arithmetic and functional
level, from the point of view of the natural bi-algebraic structure on
strata of abelian differentials. We characterise geometrically the
arithmetic points, study their distribution, and prove that in many
cases the bi-algebraic curves are the linear ones.

\end{abstract}

\maketitle
\setcounter{tocdepth}{1} 
 
\tableofcontents

\section{Introduction} \label{intro}

\subsection{Abelian differentials, their periods, and their
  bi-algebraic geometry} \label{background}
An abelian differential (or translation surface) is a pair $(C, \omega)$, where $C$
denotes a smooth irreducible complex projective curve and $\omega \in
H^0(C, \Omega^1_C) \setminus \{0\}$ is a non-zero algebraic one-form
on $C$. Its periods are the complex numbers $\int_\gamma\omega$, for 
$\gamma$ an element in the relative homology group $H_1(C^\an, Z(\omega); \ZZ)$, where $Z(\omega)$ denotes
the finite set of zeroes of $\omega$ in the compact Riemann surface $C^\an$
associated to $C$. Such a period is said to be
{\em pure} if $\gamma$ belongs to the subgroup $H_1(C^\an, \ZZ)$ of
$H_1(C^\an, Z(\omega); \ZZ)$. The goal of this paper is to study the transcendence theory of periods
of abelian differentials, both at the arithmetic and functional
level.

\medskip
At the arithmetic level: given an abelian differential $(C, \omega)$
defined over $\oQ$, we want to study the transcendence properties
of its periods. This is a
classical topic. A famous result of Schneider 
\cite{Schneider}, generalizing Siegel \cite{Siegel}, says that at
least one pure period of any abelian differential defined over $\oQ$
(meaning that both $C$ and $\omega$ are defined over $\oQ$) 
is a transcendental number. From this point of view, it is natural to consider
abelian differentials, as well as their periods, only \emph{up to scaling}: pairs $(C, [\omega])$ with
$[\omega] \in \proj H^0(C, \Omega^1_C)$ and their period lines $\left
  [\int_{\gamma_{0}} \omega, \dots, \int_{\gamma_d} \omega 
\right ]  \in \proj^{d}(\CC)$, where $d= \dim H_1(C^\an, Z(\omega);
\QQ) -1$. Our main arithmetic objects of interest are the ``least transcendent''
differentials: the $(C, [\omega])$s defined
over $\oQ$ whose period lines $\left [\int_{\gamma_{0}} \omega, \dots, \int_{\gamma_d} \omega
\right ]$ belong to $\proj^{d}(\oQ)$. Such an
abelian differential will be said to be \emph{arithmetic}.

\medskip
At the functional level: given a {\em family} of abelian
differentials, we want to study the algebraic
relations satisfied by their periods. To do so, we regroup abelian
differentials according to their combinatorial type. 
Let $g$ be a positive integer and let $\alpha$ be a partition of
$2g-2$, of length $n_\alpha$. If we ignore orbifold phenomena (see \Cref{subsec_strata}), the
abelian differentials $(C, \omega)$ with zeroes of multiplicity
$\alpha$ are parametrized by 
a smooth complex quasi-projective algebraic variety 
$\nH_\alpha$ (not necessarily irreducible), called a {\it stratum of
  abelian differentials}. It is naturally realized as a locally closed
algebraic subvariety defined over $\QQ$ of the Hodge bundle $\Omega^1
\nM_{g}$ over the coarse moduli space $\nM_g$ of smooth complex projective
curves of genus $g$. Our main geometric object of study will be the projectivization $\nS_\alpha \subset \proj \Omega^1
\nM_{g}$ of $\nH_\alpha$. This is a smooth quasi-projective variety of
dimension $d_\alpha:= 2g-2 + 
n_\alpha$, parametrizing pairs $(C,
[\omega])$ of type $\alpha$. The variety $\nH_\alpha$ is a principal
$\G_m$-bundle over $\nS_\alpha$.

\medskip
According to a remarkable theorem of Veech 
\cite[Theor. 7.15]{Ve90}, the local geometry of the complex 
manifolds $\nH_\alpha^\an$ and $\nS_\alpha^\an$, analytifications of $\nH_\alpha$
and $\nS_\alpha$ respectively, can be completely 
described in terms of periods. Let $x_0:= (C_0, \omega_0) \in
\nH_\alpha^\an$, let $\nV_{\alpha, \ZZ} := H^1(C_0, Z(\omega_0); \ZZ)$
and $\nV_\alpha=  \nV_{\alpha, \ZZ}  \otimes_\ZZ \CC$. On a small open simply connected
neighborhood $U \subset \nH_\alpha^\an$ around $x_0$, the map $D_U: U \to \nV_\alpha \simeq \CC^{d_{\alpha}+1}$
which associates to $(C, \omega) \in U$ the parallel transport in
$\nV_\alpha$ of the Betti cohomology class defined by $\omega$ in $H^1(C^\an,
Z(\omega);\CC)$, identified with the vector of periods of $\omega$, is a
bi-holomorphism. These {\em period charts} define
an {\em integral linear structure} on $\nH_\alpha^\an$, namely an atlas of 
charts with value in the complex vector space $\nV_\alpha$, whose transition
functions are locally constant elements of the 
integral linear group $\GL(\nV_{\alpha, \ZZ})$.
This integral linear structure on $\nH_\alpha^\an$ induces a {\em linear
  projective structure} on $\nS_\alpha^\an$, namely an atlas of 
charts with value in the projective space $\proj \nV_\alpha$, whose transition
functions are locally constant elements of the 
integral projective linear group $\PGL(\nV_{\alpha, \ZZ})$. These
period charts of $\nS_\alpha^\an$ are highly transcendental with
respect to the algebraic   
structure of $\nS_\alpha$: the maps $D_U$ are defined via the non-algebraic
operations of parallel transport and integration. We would like to
understand the algebraic subvarieties of $\nS_{\alpha}$ whose
periods also satisfy many algebraic relations.

\medskip
To study both the functional and arithmetic transcendence properties
of $\nS_\alpha$, we introduce the following subvarieties, which should encode most of its 
interesting geometric and arithmetic information:

\begin{defi} \label{bi-algebraic}
  \begin{enumerate}
\item 
A bi-algebraic subvariety $W \subset \nS_{\alpha}$ is an irreducible
closed algebraic subvariety $W$
of $\nS_{\alpha}$ such that $W^\an$ is algebraic in the period charts: the relative periods of $\omega$ satisfy (up
to scaling) exactly $\codim_{\nS_{\alpha}}W$ 
independent algebraic relations over $\CC$ when $(C, [\omega])$ ranges through
$W$.
\item A $\oQ$-bi-algebraic subvariety of $\nS_\alpha$ is a bi-algebraic
  subvariety $W \subset \nS_{\alpha}$ such that both $W^\an$
  (in the period charts) and $W$ are defined over $\oQ$.
  \end{enumerate}
\end{defi}

\begin{rem} \label{arithmetic points}
In particular, the $\oQ$-bi-algebraic points of $\nS_{\alpha}$ are by
definition the
arithmetic differentials. For simplicity, from now we will call \emph{arithmetic
  points} the
$\oQ$-bi-algebraic points of $\nS_\alpha$.
\end{rem}

\subsection{Linear subvarieties}

A priori, the simplest bi-algebraic subvarieties of $\nS_\alpha$ are
the {\em linear ones}: the ones for which the algebraic relations
between their periods are linear.

\begin{defi} \label{linear}
A linear subvariety $W \subset \nS_{\alpha}$ is an irreducible closed algebraic
subvariety such that $W^\an$ is (projectively) linear in the period charts.
It is $\oQ$-linear if moreover both $W$ and $W^\an$ (in the period
charts) are defined over $\oQ$.
\end{defi}

\subsubsection{Invariant linear subvarieties} \label{invariant linear}
A particular class of linear subvarieties of $\nS_\alpha$ has been studied in great
depth by dynamicists in the last twenty years: the so-called {\em
  invariant} ones. To explain the terminology, notice that the integral linear structure on
$\nH_\alpha^\an$ endows it with a a real analytic, non-algebraic
action of $\GL^+(2, \RR)$. Indeed, the group $\GL^+(2, 
\RR)$ acts naturally on $\nV_{\alpha, \ZZ}\otimes_\ZZ \CC$, by identifying $\CC$ with $\RR^2$ and extending the
natural action of $\GL^+(2, \RR)$ on $\RR^2$ coordinate-wise on
$\nV_{\alpha, \ZZ}\otimes_\ZZ \RR^2$. This action in period charts commutes with the
one of $\GL(\nV_{\alpha, \ZZ})$, hence descends to
$\nH_\alpha^\an$ (but not to $\nS_\alpha^\an$). A closed irreducible
algebraic subvariety $W$ of
$\nS_\alpha$ is said to be {\em invariant} if the analytification of its preimage in
$\nH_\alpha$ is $\GL^+(2, \RR)$-invariant. One easily shows that
a closed irreducible subvariety $W \subset \nS_\alpha$ is invariant if and
only if $W$ is linear and $W^\an$ is defined over $\RR$ (in the period
charts). Prominent examples of invariant linear subvarieties of
$\nS_\alpha$ are the famous {\em Teichm\"uller curves}: the projection
in $\nS_\alpha$ of closed $\GL^+(2, \RR)$-orbits in
$\nH_\alpha^\an$. The Teichm\"uller curves are moreover
  $\oQ$-linear. Abelian differentials $(C, [\omega])$ belonging to a
  Teichm\"uller curve are called \emph{Veech surfaces}. We refer, for
  instance, to \cite{McM03},
\cite{Calta}, \cite{Moe06}, \cite{McM07}, \cite{WrSurvey}, \cite{Moe18} for more details on Veech
surfaces and
Teichm\"uller curves. Generalizing 
Ratner's theory to this non-homogeneous setting, dynamicists proved
the following beautiful result, see \cite{McM07}, \cite{EM18}, \cite{EMM}, \cite{Wr}, \cite{F}:

\begin{theor} [McMullen in genus $2$; Eskin-Mirzakhani-Mohammadi; Wright; Filip]\label{thm_EMM;W;F}
 The topological closure of any $\GL^+(2, \RR)$-orbit in 
$\nH_{\alpha}^\an$ is the cone over an invariant $\ol{\QQ}$-linear
algebraic subvariety of $\nS_\alpha$. 
\end{theor}

\subsubsection{Linear Hurwitz spaces}\label{subsec_Hurwitz}
Let us describe nice examples of linear subvarieties of $\nS_\alpha$ which are not
invariant: the linear Hurwitz spaces. The complex manifold
$\nS_\alpha^\an$ is naturally endowed with a codimension $2g$
foliation: the  {\em isoperiodic foliation},
see~\Cref{subsec_isoper}. The isoperiodic leaf of
a point $(C_0,[\omega_0]) \in \nS_\alpha(\CC)$ consists locally of the
nearby points $(C,[\omega])$ whose vectors of \emph{pure} periods coincide with the
 one of $(C_0,[\omega_0])$. The dynamics of the isoperiodic foliation
 on $\nS_\alpha^\an$, in particular its algebraic leaves, has been recently studied in depth, see
 \cite{CDF15}. Let $E$ be an elliptic curve and let $H_{g,d}(E)$ be
 the coarse Hurwitz moduli
space classifying smooth projective curves $C$ of genus $g$ together
with a degree $d$ branched cover $C\to E$. The pullback, along the branched
cover, of the 1-form on $E$ provides a 1-form on $C$, hence a morphism
$H_{g,d}(E)\to \proj\Omega^1\nM_g$. It is not hard to see that it is
an immersion and that its image is contained in an isoperiodic
leaf. One easily checks that the connected
components of $(H_{g,d}(E)\cap\nS_\alpha)^\an$ are saturated with
respect to the foliation, hence coincide with closed isoperiodic
leaves, and thus are linear subvarieties. It is proved in
\cite{CDF15} that for $g>2$ these are the only algebraic leaves of the
isoperiodic foliation. If the elliptic curve $E$ is not
 defined over $\Qbar$, these linear subvarieties are not defined over
 $\Qbar$. In particular they are not invariant.

\subsubsection{Other examples}

In \cite[Section 6]{Moe08}, M\"oller exhibits finitely many
non-isoperiodic linear subvarieties not defined over $\RR$: they are
constructed from the families of cyclic covers of $\proj^1$ studied by
Deligne-Mostow \cite{DM86}. 


\section{Results and conjectures}
The bi-algebraic format for $\nS_\alpha$ introduced in \Cref{bi-algebraic} is a special
instance of the general 
bi-algebraic format described for instance in \cite{KUY2}. Given a complex algebraic variety $S$ with an infinite
topological fundamental group, this format proposes to emulate an
algebraic structure on the universal cover 
$\widetilde{S^\an}$ of $S^\an$ using periods of algebraic 
differential forms on $S$; to analyse the bi-algebraic and
$\oQ$-bi-algebraic subvarieties of $S$;  to study the trancendence  
properties of the uniformizing map $\pi: \widetilde{S^\an} \to
S^\an$ with respect to the emulated algebraic structure on
$\widetilde{S^\an}$ and the algebraic structure of $S$ (Ax-Lindemann and Ax-Schanuel heuristic); and to
analyse the distribution of the arithmetic (= $\oQ$-bi-algebraic) points (Zilber-Pink
heuristic). This format has first been successfully applied to tori, abelian varieties
and Shimura varieties, see \cite{KUY2} for a survey and references;
then to general varieties $S$ endowed with a variation of (possibly mixed) Hodge structures \cite{Klingler},
\cite{BT19}, \cite{BKT}, \cite{KO}, \cite{BKU}. In this paper we go
one step further, moving to an even less homogeneous context which creates considerable new difficulties.

\subsection{Arithmetic points}

\subsubsection{Geometric characterization of the arithmetic points}
Our first main result in this paper is the geometric elucidation
of the arithmetic points in $\nS_\alpha$. We will use the following

\begin{defi} \label{degree}
  Let $(C, [\omega])$ be an abelian differential. We define:
\begin{enumerate}
  \item the abelian variety $A_{[\omega]}$, as the smallest factor of the Albanese $\Alb(C)$ (canonically
    identified with the Jacobian $\Jac(C))$ whose tangent bundle contains
    $\omega$. The \emph{degree} $d_{[\omega]}$ of $(C, [\omega])$ is the dimension
    $\dim_\CC A_{[\omega]}$. Thus $1 \leq d_{[\omega]} \leq
g(C)$, where $g(C)$ is the genus of $C$.
   \item the line $[\omega_{A_{[\omega]}}]   \in \proj H^0(A_{[\omega]},
    \Omega^1_{A_{[\omega]}})$, as the unique
    line such that $[\omega] = [\alb_{[\omega]}^* \,
    \omega_{A_{[\omega]}}]$. Here $\alb_{[\omega]}: C \to
      A_{[\omega]}$ is the composition of the Albanese map $\alb: C \to \Alb(C)$ with the
  projection of $\Alb(C)$ onto $A_{[\omega]}$ (the maps $\alb$ and $\alb_{[\omega]}$ are
  uniquely defined up to a translation).
  \end{enumerate}
  \end{defi}

  \noindent We also refer to \Cref{CM} for our (standard) terminology concerning complex
  multiplication (CM).

\begin{theor} \label{arithmetic points:geometry}
A point $(C, [\omega]) \in \nS_{\alpha}(\oQ)$ is arithmetic if and only
if the following conditions are satisfied:

\begin{enumerate}
\item The complex curve $C$ is defined over $\oQ$: $C= C_{\oQ} \otimes_{\oQ}
\CC$.

\item``The differential $[\omega]$ is an eigenform for complex multiplication''; namely:
\begin{enumerate}
  \item[(a)] The abelian variety $A_{[\omega]}$ has complex
    multiplication and is isotypic.
  \item[(b)] The line $[\omega_{A_{[\omega]}}] \in \proj
    H^0(A_{[\omega]}, \Omega^1_{A_{[\omega]}})$ is an eigenline
      for the $K$-action on $H^0(A_{[\omega]},
      \Omega^1_{A_{[\omega]}})$, where $K$ denotes the CM field center of
      $\End_\QQ A_{[\omega]}$.
  \end{enumerate}
  
\item Given any two points $x$ and $y$ in $Z([\omega])$, the
  difference $\alb_{[\omega]} (x) - alb_{[\omega]} (y)$ is a torsion
  point of $A_{{[\omega]}}$.
\end{enumerate}
\end{theor}

The main tool in the proof of \Cref{arithmetic points:geometry} is,
as for most results nowadays in arithmetic transcendence theory,
W\"ustholz' Analytic Subgroup Theorem \cite{wus}.

\subsubsection{Arithmetic points and Veech surfaces} \label{Veech}
There is a striking similarity between \Cref{arithmetic points:geometry} characterising
arithmetic points, and properties of the Veech surfaces.  In \cite[Theorem 2.7]{Moe06} and
  \cite[Theorem 3.3]{Moe06a}, M\"oller shows that Veech surfaces 
  satisfy conditions similar to, but essentially weaker than, the ones of
  \Cref{arithmetic points:geometry}: if an abelian
  differential $(C,[\omega])$ is a Veech surface, then:
  \begin{enumerate}
    \item[(1V)] The factor
  $A_{[\omega]}$ of $\Alb(C)$ has {\em real} multiplication by a totally real
  field $K_0$ satisfying $[K_0: \QQ]= d_{[\omega]}$ (in particular
  $A_{[\omega]}$ is isotypic);
  \item[(2V)] The line
  $[\omega_{A_{[\omega]}}]$ is an eigenline for the $K_0$-action on $H^0(A_{[\omega]},
  \Omega^1_{A_{[\omega]}})$;
  \item[(3V)] Condition (3) of \Cref{arithmetic points:geometry}
    holds.
  \end{enumerate}
  Moreover the degree $d_{[\omega]}$ of a Veech surface $(C, [\omega])$ coincides
  with the degree of the Teichm\"uller curve $\Gamma \backslash \fH
  \subset \nS_\alpha^\an$ it generates, 
    defined as the degree of its trace field $\QQ[\tr \, \gamma, \;
    \gamma \in \Gamma]$.

    \medskip
It follows immediately from \Cref{arithmetic points:geometry} that any
    arithmetic point $(C, [\omega])$ with $A_{[\omega]}$ simple also
    satisfy the conditions $(1V)$, $(2V)$ and $(3V)$ (for $K_0$ the maximal
    totally real subfield of the CM field $K$). On the other hand,
    these conditions are in general not sufficient for an abelian
  differential $(C, \omega)$ to be a Veech surface. The following
  proposition clarifies the relation between arithmetic points and
  Veech surfaces, showing that arithmetic points of degree $1$ are
  familiar, while arithmetic points of higher degree are mysterious
  and do not seem related to the $\GL^+(2, \RR)$-action on $\nH_\alpha^\an$:

  \begin{prop} \label{arithm-versus-Veech}
   The arithmetic points of $\nS_\alpha^\an$ of degree
    $1$ are Veech surfaces. On the other hand there exist arithmetic points of
  degree at least $2$ which are not Veech surfaces.
  \end{prop}

  
  \subsubsection{Distribution of arithmetic points}
  
The arithmetic points are many:

\begin{prop} \label{arithmetic points:density}
Arithmetic points of degree $1$ are (analytically) dense in $\nS_{\alpha}^\an$ for all $\alpha$.
\end{prop}

\begin{rem}
  This is where considering abelian differentials {\it up to scaling} is
  crucial. If we were to consider strata $\nH_\alpha$ of abelian differentials
  without scaling,
  with their natural integral linear structure, Schneider's
  theorem implies that there are no arithmetic points
  at all in $\nH_\alpha$!
\end{rem}

On the other hand, in stark contrast with the
bi-algebraic geometry of tori, abelian varieties or Shimura varieties,
we show that the 
arithmetic points are in general not Zariski dense in the $\oQ$-bi-algebraic
  subvarieties of $\nS_\alpha$. The situation here is thus similar
  with the one for general variations of Hodge structures. Our second main result in this paper is the
  following (see \Cref{sec_counterexample} for the definition of the
  rank and the
  degree of a general invariant linear subvariety; any Teichm\'uller
  curve has rank $1$):
  
\begin{theor} \label{prop_counterexample}
Any invariant linear subvariety of $\nS_\alpha$ of rank~$1$ and degree at least~$2$
does not contain a Zariski-dense set of arithmetic points. In
particular any Teichm\"uller curve of degree at least $2$ contains
only finitely many arithmetic points.

In general, an invariant linear subvariety of $\nS_\alpha$ of
rank~$k\geq 1$ and
degree~$d >1$ does not contain a Zariski-dense set of arithmetic points
of degrees at least $kd$.
\end{theor}

\noindent
In addition to \Cref{arithmetic points:geometry} and the results of
\cite{Moe06}, \cite{Fi}, \cite{Fil17} and 
\cite{EFW18} on invariant linear varieties, the main tool in the proof
of \Cref{prop_counterexample} is the Andr\'e-Oort conjecture for
mixed Shimura varieties whose pure part is of abelian type \cite{Ts18}, \cite{Gao}.

\subsection{Bi-algebraic subvarieties}
Let us now turn to the purely geometric aspects of the bi-algebraic
geometry of $\nS_\alpha$. For tori, abelian varieties or Shimura
varieties, its bi-algebraic subvarieties are its weakly special
subvarieties: they are defined by group-theoretic conditions and may be thought as the most
``linear'' algebraic subvarieties.
We propose the following:

\begin{conj} \label{linear conj}
The bi-algebraic subvarieties of $\nS_\alpha$ are the linear ones.
\end{conj}

\noindent
In words: if the periods of an algebraic family of abelian
differentials satisfy as many algebraic relations as possible, then
these algebraic relations are linear.

\medskip
The absence of homogeneity for $\widetilde{\nS_\alpha^\an}$ makes \Cref{linear conj}
of a completely different order of difficulty than for tori, abelian
varieties or Shimura varieties. Our third main result is a
proof of \Cref{linear conj} in two cases: for bi-algebraic curves on which $\omega$
does not vary in a constant local subsystem (see the
condition~$(\star)$ of \Cref{cond_star}); and for 
bi-algebraic subvarieties of $\nS_{\alpha}$ contained
in an isoperiodic leaf.

\begin{theor} \label{bi-algebraic=linear1}
The bi-algebraic curves in $\nS_{\alpha}$ satisfying condition
  $(\star)$ are linear.
\end{theor}

\begin{theor} \label{bi-algebraic=linear2}
The bi-algebraic subvarieties of $\nS_{\alpha}$ contained in an isoperiodic leaf are linear. 
\end{theor}

\noi
The main tools for proving \Cref{bi-algebraic=linear1} and
\Cref{bi-algebraic=linear2} are a detailed analysis of the monodromy
action, using the classical results of Satake \cite{Sa65}, as well as
the Ax-Schanuel conjecture for abelian varieties \cite{Ax72}.

\medskip
Thanks to \Cref{bi-algebraic=linear1} and \Cref{bi-algebraic=linear2},
\Cref{linear conj} holds true in genus~$2$:
\begin{theor} \label{genus2}
Any bi-algebraic curve of $\nS_{1,1}$ or $\nS_{2}$ either coincide
with a linear Hurwitz space, or satisfy condition $(\star)$. In
particular it is linear.
\end{theor}

\begin{rem}
Notice that \cite[Theorem 7.1]{Moe08} fully describes linear curves in
genus 2. Thus \Cref{genus2}, in combination with this description,
provides a complete classification of the bi-algebraic
curves in genus $2$.
\end{rem}

In genus~$3$, we show in \Cref{genus3} an example of a bi-algebraic curve in $\nS_4$
which is not a linear Hurwitz space, nor satisfies condition
$(\star)$, but is still linear.

\subsection{Conjectures} \label{conj}
Let us propose some questions and
conjectures for $\nS_{\alpha}$, which are suggested
by the bi-algebraic format and which place the previous results in their proper context.

\medskip
First of all, in view of \Cref{prop_counterexample}, it would be interesting to decide the following

\begin{Question} \label{Question}
Does any $\oQ$-bi-algebraic subvariety $\nS_\alpha$
contain at least one arithmetic point?
\end{Question}

We check in \Cref{bow} that this holds true for the unique known series of primitive
  Teichm\"uller curves generated by Veech surfaces of unbounded
  genera: the triangle group series of Bouw-M\"oller \cite{BM10}. One
  also easily checks this is true
 for the $\oQ$-linear Hurwitz spaces. 

\medskip
Even if not every $\oQ$-bi-algebraic subvariety of
$\nS_\alpha$ does contain a Zariski-dense set of arithmetic points, we
still conjecture that the converse holds true:

\begin{conj} [Andr\'e-Oort for $\nS_{\alpha}$] \label{AO}
Let $S \subset \nS_{\alpha}$ be an irreducible algebraic subvariety of
$\nS_{\alpha}$ containing a Zariski-dense set of arithmetic points. Then
$S$ is $\oQ$-bi-algebraic (hence $\oQ$-linear in view of \Cref{linear conj}).
\end{conj}

Both \Cref{prop_counterexample} and \Cref{AO}, as well as the main result
of \cite{EFW18} stating that all but
finitely many linear invariant subvarieties of $\nS_\alpha$ have degree
at most $2$, are instances of a general \emph{Zilber-Pink conjecture}
for \emph{atypical intersections} 
in $\nS_\alpha$. We won't write about it in detail here and will 
come back to it in a future article. Let us just state the following,
which we see as the core of this conjecture:

\begin{conj} \label{finiteness}
  Any stratum $\nS_\alpha$ contains only finitely many arithmetic
  points of degree at least $3$.
\end{conj}

\noi
In the same way as for other bi-algebraic structures (see \cite{UY2}, \cite{PT}, \cite{KUY},
\cite{Klingler}), the main geometric step towards
the Zilber-Pink conjecture is an \emph{Ax-Schanuel conjecture} for
$\nS_\alpha$, of which we will describe here only the following particular case:
\begin{conj} [Ax-Lindemann for $\nS_{\alpha}$] \label{AL}
Let $V \subset \widetilde{{\nS}_{\alpha}^\an}$ be an irreducible algebraic
subvariety of the universal cover $\widetilde{{\nS}_{\alpha}^\an}$ of
  ${\nS}_{\alpha}^\an$, namely a closed irreducible analytic subvariety of
$\widetilde{{\nS}_{\alpha}^\an}$ which is algebraic in the period
charts. Then the Zariski-closure $\overline{\pi(V)}^\Zar$ of its
projection in $\nS_\alpha$ is
bi-algebraic  in $\nS_{\alpha}$ (hence linear in view of \Cref{linear conj}).
\end{conj}



\subsection{Organization of the paper}

\Cref{prelim} fixes the notations and complements the introduction with
the details we need later, including the mixed Hodge theory of strata and its
relation to their bi-algebraic geometry. 
In \Cref{SP:geometry} we prove the geometric characterization \Cref{arithmetic points:geometry}
of the arithmetic points. \Cref{density} discusses the distribution of the arithmetic
points and \Cref{prop_counterexample}. \Cref{linearity} proves
\Cref{bi-algebraic=linear1} and \Cref{bi-algebraic=linear2}. 

\subsection{Acknowledgments}
B.K. would like to thank Madhav Nori, for the many discussions on
this topic a long time ago.


\section{Preliminaries} \label{prelim}

\subsection{Strata}\label{subsec_strata}

Let $(C, [\omega])$ be an abelian differential (up to scaling) of
genus $g \geq 1$. The divisor $\Div([\omega])$ on $C$ of zeroes of
$[\omega]$ can be uniquely written $\sum_{i=1}^n
\alpha_i x_i $, with $\alpha_i \in \NN^*$ and the $x_i$s are pairwise
distinct in $C$.  The
set of integers $\alpha\coloneqq \{\alpha_1, \dots, \alpha_n\}$ is called the type
of $(C, [\omega])$. As the weight $|\alpha |:= \sum_{i=1}^n \alpha_i$ of 
$\alpha$ coincides with the degree of $\Div(\omega)$, it satisfies $|\alpha| =
2g-2$.

\begin{rem}
The algebro-geometric notion of abelian differential is equivalent to the differential
geometric notion of {\em translation surface},  i.e. a compact Riemann
surface $\Sigma$, a finite union of points $Z 
\subset \Sigma$, and an atlas of charts for $\Sigma \setminus Z$ with value in
$\CC$ whose transition functions are locally constant translations (a so-called
{\em translation structure} on $\Sigma \setminus Z$) such that the cone
angle at each point of $\Sigma$ is a positive integral multiple of $2\pi$. The
equivalence is obtained by associating to $(C, \omega)$ the pair
$(\Sigma:= C^\an, Z:= Z(\omega))$, the translation structure
being given by locally integrating $\omega$ on $C^\an \setminus
Z(\omega)$. This differential geometric point of view will play no role in
this paper.
\end{rem}

Let $\nM_g$ denote the coarse moduli space of
smooth projective curves of genus $g$. The coarse moduli space of
non-zero abelian differentials (up to scaling) of genus $g$ is the 
projectivization $\proj \Omega^1 \nM_g$ of the Hodge bundle $\Omega^1
\nM_g$ on $\nM_g$ (whose fiber at a closed point $C \in \nM_g$
is the space of algebraic
one-forms on $C$).  The algebraic variety $\proj
\Omega^1 \nM_g$ is defined over $\QQ$ and has dimension $4g-4$. It
is naturally stratified according to the type of the one-forms:
$$ \proj \Omega^1 \nM_g= \coprod_{\substack{\alpha = \{\alpha_1,
    \dots, \alpha_{n}\} \\ \sum_i \alpha_{i} = 2g-2}} \nS_{\alpha}$$
where the stratum $\nS_{\alpha}$ parametrizing abelian differentials
of type $\alpha$ is defined set-theoretically as the abelian
differentials $(C, [\omega]) \in \proj \Omega^1 \nM_g$ for which there
are pairwise distinct points $x_i \in C  \textnormal {\, with\;} \Div([\omega]) =
\sum_{i=1}^n \alpha_i x_i$.
This is a quasi-projective variety defined over $\QQ$ and has dimension
$2g+n-2$. It is not connected in
general, see \cite{KZ03}. 
All the coarse moduli spaces introduced above have orbifold singularities. 
To define the period coordinates around every point we will need to
work with fine moduli spaces: this can be achieved by passing to a
finite ramified cover, for instance introducing a level-$\ell$
structure, with $\ell\ge 3$, on the curves underlying the abelian
differentials.  
{\em For the rest of the text, we fix a level $\ell\ge3$ and, by abuse of notations, we call stratum
  and denote by $\nS_\alpha$ any connected component of the level-$\ell$ ramified cover of the
  $\nS_\alpha$ defined above}.  
With this convention, $\nS_\alpha$ is a smooth, connected,
quasi-projective variety defined over $\Qbar$. In the same way, we
suppress the integer $\ell$ from 
  the notation of the moduli space of level-$\ell$ curves: we denote
  by $\cM_g$ (resp. $\cM_{g,n}$) the fine moduli space of level-$\ell$
  smooth projective curves of genus $g$ (resp. with $n$ distinct
  marked points). 
Our results do not depend on the choice of the level.

\medskip
Let us mention two variants of the above definitions which we will use.
First, it will be sometimes convenient to work with the moduli space
$\nH_\alpha$ of abelian differentials {\it without scaling}, with
zeros of type $\alpha$ (a $\Gm$-bundle over $\nS_\alpha$). It will
also be useful to consider a version of $\nS_\alpha$ where the zeroes
of the forms are marked. Let $\cM_{g, n}$ be the fine moduli space of level-$\ell$ smooth projective curves of
genus $g$ with $n$ distinct marked points. The moduli space
of canonical divisors of type $\alpha$ is the locally closed
subvariety of $\cM_{g, n}$ defined set theoretically by
\begin{equation*} \label{moduli of canonical divisors}
  \cS_\alpha:= \left \{ [C, x_1, \dots, x_n] \in \cM_{g, n} \; | \;
      \cO_C(\sum_{i=1}^n \alpha_i \, x_i) \simeq \Omega^1_C \right
    \} \stackrel{\iota}{\hookrightarrow} \cM_{g, n}
    \;\;.
  \end{equation*}
  
Notice that if $\alpha = (m_1, \dots, m_1, m_2, \dots, m_2, \dots,
m_r, \dots, m_r)$ with $m_i \neq m_j$ for $i \neq j$, then the
product of symmetric groups $\Sigma_{n_1}\times \dots \times
\Sigma_{n_{r}}$ acts freely on $\cS_\alpha$; and $\cS_\alpha$ is a
finite \'etale cover of $\nS_\alpha$.

\subsection{Hodge theory of $n$-pointed curves} 
We first recall the classical mixed Hodge theory of curves with
marked points. Let $C$ be a smooth complex projective curve of
genus $g$. The Betti cohomology group $H^1(C, \ZZ)$ is a pure $\ZZ$-Hodge
structure of weight $1$. If $Z=\{x_1, \dots, x_n\}$ is a set of $n$ distinct
points on $C$, the long exact sequence of relative cohomology
for the pair $(C, Z)$
\begin{multline*} 
\dots \to H^0(C, \ZZ) \to H^0(Z, \ZZ) 
\to H^1(C, Z, \ZZ) \to H^1(C, \ZZ) \to
H^1(Z, \ZZ)=0 
\end{multline*}
defines a short exact sequence in the category of $\ZZ$-mixed Hodge
structures (abbreviated $\ZZ$MHS)
\begin{equation} \label{MHS}
0 \to \tilde{H}^0(Z, \ZZ) \simeq \ZZ(0)^{n-1} \to H^1(C, Z, \ZZ) \to
H^1(C, \ZZ) \to 0\;\;,
\end{equation}
where $\tilde{H}^\bullet$ denotes the reduced cohomology.
The $\ZZ$-mixed Hodge structure on $H^1(C, Z, \ZZ)$ is of
type $(0,0)$, $(1,0)$ and $(0,1)$:

- $\Gr^W_0 H^1(C,
Z; \QQ) = \QQ(0)^{n-1}$ and $\Gr^W_1 H^1(C,
Z; \QQ) = H^1(C, \ZZ)$.

- the only non-trivial piece of the Hodge filtration is $F^1H^1(C, Z;
\CC)$, defined as the subspace $H^0(C, \Omega^1_C)$ of $H^1(C, Z;\CC)$.

\begin{rem} \label{dual1}
As $C$ is a smooth projective curve, Poincar\'e duality provides various
identifications of Hodge structures $H^1(C, Z, \ZZ) \simeq H_1(C-Z; \ZZ)(-1) \simeq H^1_c(C
-Z, \ZZ)$, while the dual Hodge structure $H^1(C, Z, \ZZ)^\vee$
identifies with $H^1(C-Z, \ZZ)(1) \simeq H_1(C,
Z,\ZZ)$. As a result, any natural $\ZZ$-mixed Hodge structure
associated with the pair $(C, Z)$ coincides, up to twist and
duality, with (\ref{MHS}). For example the $\ZZ$-mixed Hodge structure
$$ 0 \to H^1(C, \ZZ)(1) \to H^1(C-Z,\ZZ)(1) \to \tilde{H}^0(Z,\ZZ) \to
0$$
is dual to (\ref{MHS}).
The only non-trivial piece of the Hodge filtration on $H^1(C-Z,\CC)$ is $F^1H^1(C-Z,
\CC)$, defined as the space $H^0(C, \Omega^1_C(\log Z))$ of
logarithmic one-forms on $C$ with poles at $Z$, which is the natural
annihilator of $F^{0}H^1(C, Z;
\CC)(1)\simeq H^0(C, \Omega^1_C)$.
\end{rem}


\subsection{The variation of $\ZZ$-mixed Hodge structure $\bV_\alpha$
  on $\nS_\alpha$}

\subsubsection{}

When $C$ varies through $\cM_g$ the cohomology $H^1(C, \ZZ)$ defines a
$\ZZ$-variation of Hodge structure ($\ZZ$VHS) $\bV_{g, \ZZ}$ on $\cM_g$ of
weight one. If $f: \mathcal{C}_{g} \to \cM_{g}$ denotes the universal
smooth projective curve of genus $g$ and level-$\ell$ then $\bV_{g, \ZZ}= R^1f_* \ZZ$.
Varying $(C, [\omega])$ through $\nS_\alpha$, one
obtains similarly that $\bV_{\alpha, \ZZ}$ is an admissible, graded polarized, variation of $\ZZ$-mixed Hodge
structure ($\ZZ$VMHS). The weight filtration has two steps:
\begin{equation} \label{vmhs}
  0 \to W_0 \bV_{\alpha, \ZZ} \simeq \ZZ(0)^{n-1} \to \bV_{\alpha, \ZZ} \to
  \Gr_1^W \bV_{\alpha, \ZZ} = p^* \bV_{g, \ZZ} \to 0\;\;,
  \end{equation}
where $p: \nS_\alpha \to \cM_g$ denotes the canonical projection.

\begin{rem}
  If one replaces $\nS_\alpha$ by its finite \'etale cover
  $\cS_\alpha \stackrel{\iota}{\hookrightarrow} \cM_{g,n}$, then the
  $\ZZ$VMHS $\bV_{\alpha, \ZZ}$ on $\cS_\alpha$ coincides with
  $\iota^{-1} R^1{f^0}_!\ZZ_{\mathcal{C}_{g,n}^0}$, where $f^0:
\mathcal{C}_{g,n}^0 \to \cM_{g,n}$ is the open curve complement in $f:
\mathcal{C}_{g,n} \to \cM_{g,n}$ (the universal smooth projective 
$n$-pointed curve of genus $g$ and level-$\ell$) of the
canonical sections $x_i: \cM_{g,n} \to \mathcal{C}_{g,n}$, $1 \leq i
\leq n$.
\end{rem}

\begin{rem}
It is interesting to notice that the variety $\cS_\alpha$ defined set-theoretically by (\ref{moduli of canonical
  divisors}) can itself be defined in a purely Hodge theoretic way. 
Let $\Pic^0_{\mathcal{C}_g/\cM_g} \to \cM_g$ denotes the relative Picard
scheme whose fiber at $C \in \cM_g$ is the abelian variety $\Pic^0(C)$
parametrizing degree zero line bundles on $C$, and $p: \cM_{g, n} \to
\cM_g$ the natural map forgetting the marking. The abelian scheme
$p^*\Pic^0_{\mathcal{C}_g/\cM_g} \to \cM_{g, n}$ has two natural
sections: the identity section $e$ and the section $s_\alpha$ defined by associating to $(C, x_1, \dots,
  x_{n}) \in \cM_{g,n}$ the degree zero line bundle $\Omega^1_C(-\sum_{i=1}^{n}
\alpha_i x_i)$ on $C$. The variety $\cS_\alpha$ is the subvariety of
$\cM_{g,n}$ defined by the Cartesian diagram

$$
\xymatrix{
\cS_\alpha \ar[d] \ar[r] &  \cM_{g,n} \ar[d]^<<<<e \\
\cM_{g,n} \ar[r]_>>>>>>{s_{\alpha}}  & p^* \Pic^0_{\mathcal{C}_g/\cM_g} \;. } 
$$
The variety $\cM_{g,n}$ has dimension $3g-3+n$, thus $p^*
\Pic^0_{\mathcal{C}_g/\cM_g}$ has dimension $4g-3+n$. The two sections
$s_\alpha(\cM_{g,n})$ and $e(\cM_{g,n})$ are thus of codimension $g$. On the other
hand their intersection $\cS_{\alpha}$ has dimension $2g-2+n$ hence is
of codimension $2g-1$ in $p^*\Pic^0_{\mathcal{C}_g/\cM_g}$. This shows
that $s_\alpha(\cM_{g,n})$ and $e(\cM_{g,n})$ are not transverse.
Over $\CC$ the Picard variety $\Pic^0(C)$ is canonically
isomorphic to the group of extensions $\Ext^1_{\ZZ\MHS}( H^1(C;\ZZ), \ZZ(0))$ in
the abelian category of $\ZZ$-mixed Hodge structures.
The section $s_\alpha: \cM_{g,n} \to p^*\Pic^0_{\mathcal{C}_g/\cM_g}$ can be thought as an
element of the group of extensions  $\Ext^1_{\ZZ\VMHS^{\textnormal{adm}}_{\cM_{g,n}}} (p^* (R^1f_*\ZZ),
\ZZ(0))$ in the abelian category of admissible variations of $\ZZ$MHS
on $\cM_{g,n}$
i.e. a normal function; and $\cS_\alpha$ as the zero locus of
this normal function.
\end{rem}


\subsection{The period map \texorpdfstring{$\Phi: \nS_{\alpha} \to
    \fA_g^{(n-1)} $}{Phi}} \label{HodgeDiagram}

The classifying space for weight one $\ZZ$VHS of dimension $2g$ is the
Shimura variety $\cA_g$ moduli space of principally polarized abelian
varieties of dimension $g$. The classifying space for the $\ZZ$VMHS
extensions of a weight one $\ZZ$VHS  of dimension $2g$ by $\ZZ(0)$ is
the mixed Shimura variety $\mathfrak{A}_g$, the universal principally
polarized abelian variety of dimension $g$ over $\cA_g$.
Hence the classifying space for the $\ZZ$VMHS extensions of a weight
one $\ZZ$VHS  of dimension $2g$ by $\ZZ(0)^{n-1}$ 
is the mixed Shimura variety $\mathfrak{A}_g^{(n-1)}\coloneqq \mathfrak{A}_g \times_{\cA_g}
\mathfrak{A}_g  \times_{\cA_g} \dots \times_{\cA_g}
\mathfrak{A}_g$ (product of $(n-1)$ factors). The $\ZZ$VMHS $\bV_{\alpha, \ZZ}$ on $\nS_\alpha$ is
classified by the period map
\begin{equation} \label{defi period map}
  \Phi: \nS_\alpha \to
  \mathfrak{A}_g^{(n-1)}\;\;.
\end{equation}

By a refined
version of the classical Abel-Jacobi theorem (see
\cite[Theor. 7.2]{ArOh97}), the period map is quasi-finite.

\subsection{Integral projective structure on $\nS_\alpha$} \label{linear structure}

Recall that an {\em integral projective structure} on a complex manifold $M$ can be defined equivalently
\begin{itemize}
  \item[(i)] as a maximal atlas of charts for $M$ with value in $\proj
    (\nV_\ZZ \otimes_\ZZ \CC)$ whose transition functions are locally
    constant elements of the integral projective linear group
    $\PGL(\nV_\ZZ)$; here $V_\ZZ$ denotes a finite free $\ZZ$-module
    of rank the (complex) dimension of $M$.
\item[(ii)] as a local
biholomorphism
$
D: \tilde{M} \to \proj(\nV_\ZZ \otimes_\ZZ \CC)
$,
called the developing map of the projective structure,
which is equivariant under a monodromy representation $\rho:
\pi_1(M, x_0) \to \PGL(\nV_\ZZ)$. Here $\pi:
\widetilde{M} \to M$ is the universal 
cover of $M$ 
at $x_0$, and $D$ is obtained by glueing the local projective charts.
\end{itemize}

In the case of $\nS_\alpha^\an$, let $\bV_{\alpha, \ZZ}$ be the $\ZZ$-local system on $\nS_\alpha^\an$ whose
fiber at any point $(C, [\omega]) \in \nS_\alpha^\an$ is the relative cohomology
group $H^1(C, Z([\omega]); \ZZ)$. Let $\bV_\alpha$ be its
complexification and $(\cV_\alpha^\an, \nabla^\an)$ the associated
complex analytic integrable connection on $\nS_\alpha^\an$. Given
$U\subset \nS_\alpha^\an$ a simply
connected neighbourhood of a point
$x_0:= (C_0,[\omega_0])$ in $\nS_{\alpha}^\an$, the local
  period chart of the introduction is the local biholomorphism
$$
  D_U: U \stackrel{[\underline{\omega}]}{\rightarrow} {\proj \cV_\alpha^\an}_{|U}
  \stackrel{\varphi}{\rightarrow} \PV_{\alpha, x_0}:= \proj H^1(C_0, Z([\omega_0]); \CC)
  \stackrel{\gamma}{\simeq} \proj^{2g+n-2}\CC $$ which to a point $x \in U$
  associates the vector of periods $$\left[ \int_{\gamma_1}
\varphi(\omega_x), \dots,  \int_{\gamma_{2g+n-1}}
\varphi(\omega_x) \right]\in \proj^{2g+n-2}\CC\;\;.$$ Here $[\underline{\omega}]$ is the tautological
section of $\proj \cV_\alpha$ which to a point $(C, [\omega]) \in \nS_{\alpha}^\an$
associates the class of $[\omega]_{\Betti}$ in $\proj H^1(C, Z([\omega]), \CC)$; the map
$\varphi$ is the parallel transport with respect to $\nabla^\an$ of
${\proj \cV_\alpha^\an}_{|U}$ on its central fiber $ \PV_{\alpha, x_{0}}$; and
$\gamma$ is the identification of $\PV_{\alpha, x_{0}}$ with
$\proj^{2g+n-2}\CC$ provided by the choice of an integral basis $(\gamma_i)_{1 \leq i \leq
2g+n-1}$ of $H_1(C_0, Z([\omega_0]); \ZZ)$. Changing
the base point $(C_0, \omega_0)$ results in a locally constant change
of coordinates with integral coefficients  between the two
corresponding charts.

\medskip
More globally, if $\pi: \widetilde{\nS _\alpha^\an} \to
\nS_\alpha^\an$ denotes the universal cover of $\nS_\alpha^\an$,  the
canonical trivialization of the local system $\pi^*  
\bV_{\alpha, \ZZ}$ on the simply connected space 
$\widetilde{\nS_\alpha^\an}$ induces a decomposition of the associated
vector bundle  $\widetilde{\cV_\alpha^\an}$ into a product $\widetilde{\nS_\alpha^\an} \times
\nV_{\alpha}$, where $\nV_\alpha := \nV_{\alpha, \ZZ} \otimes_\ZZ \CC$
and $\nV_{\alpha, \ZZ}:= H^0(\wti{\nS^\an_\alpha}, \pi^* \bV_{\alpha, \ZZ})$ is isomorphic to any of
the fibers $\nV_{\alpha, \ZZ, x_{0}}$. The developing map
\begin{equation} \label{eq}
D: \widetilde{\nS_\alpha^\an} \to \PV_\alpha
\end{equation}
of the projective structure on $\nS_\alpha$ is obtained as the composition of the section
$\widetilde{[\underline{\omega}]}: \widetilde{\nS_\alpha^\an} \to
\widetilde{\proj \cV^\an_\alpha}$ lifting
$[\underline{\omega}]$, followed by the parallel transport projection
$\widetilde{\proj \cV_\alpha^\an}\simeq
\widetilde{\nS_\alpha^\an}  \times \PV_\alpha \to \PV_\alpha$. It is naturally
equivariant under $\pi_1(\nS_\alpha^\an)$, we denote by $\rho:
\pi_1(\nS_\alpha^\an) \to \PGL(\nV_{\alpha, \ZZ})$ the associated monodromy.


\subsection{The bi-algebraic geometry of
  \texorpdfstring{$\nS_\alpha$}{PH}}

Recall from \cite[Section 4]{KUY} that:

\begin{defi} \label{bialg}
A bi-algebraic structure on a connected complex algebraic variety $S$
is a pair 
$$ (D: \wti{S^\an} \to X, \quad \rho: \pi_1(S^\an) \to \Aut(X))\;\;$$ 
where $\wti{S^\an}$ denotes the topological universal cover of $S$, $X$ is a complex
algebraic variety, $\Aut(X)$ its group of algebraic
automorphisms, $\rho: \pi_1(S^\an) \to \Aut(X)$ is a group morphism
(called the holonomy representation) and 
$D$ is a $\rho$-equivariant holomorphic map (called the developing map).
\end{defi}

\begin{defi} \label{bialgebraic2}
A $\oQ$-bi-algebraic structure on a complex algebraic variety
$S$ is a complex bi-algebraic structure $(D: \wti{S^\an} \to X, \rho:
\pi_1(S^\an) \to \Aut(X))$ such that:

(1) $S$ is defined over $\overline{\QQ}$.

(2) $X = X_{\oQ} \otimes_{\oQ} \CC$ is defined over $\overline{\QQ}$ and the
  homomorphism $\rho$ takes value in $\Aut_{\oQ} X_{\oQ}$.
\end{defi}

\begin{defi}
The {\em linear $\oQ$-bi-algebraic structure} on $\nS_\alpha$ is the
one defined by $(D, \rho)$ as in (\ref{eq}).
\end{defi}

Following the general
format of \cite[Section 4]{KUY}, we define:

\begin{defi} \label{defi1}
A complex analytic subvariety $Y \subset
\widetilde{{\nS}_{\alpha}^\an}$ is said to be algebraic if it is
an irreducible complex analytic component of $D^{-1}(Z)$, for $Z \subset
\PV_{\alpha}$ a complex algebraic subvariety. It is said to be defined over $\oQ$ if $Z$ is.
\end{defi} 

\begin{rem} \label{algebraic model}
Notice that in \Cref{defi1} one can equivalently replace $Z$ by the
Zariski-closure $\overline{D(Y)}^\Zar$ of $D(Y)$ in the projective
space $\PV_\alpha$. We call $\overline{D(Y)}^\Zar$ {\em the
algebraic model} of $Y$.
\end{rem}

\begin{defi} \label{bi-algebraic1}
A bi-algebraic subvariety $W \subset \nS_{\alpha}$ is an irreducible algebraic subvariety $W$
of $\nS_{\alpha}$ such that $W^\an$ is the projection $\pi(Y)$ of an 
algebraic subvariety $Y$ of $\widetilde{{\nS}_{\alpha}^\an}$ (in the sense
of \Cref{defi1}). It is $\oQ$-bi-algebraic if moreover $W$ and $Y$ are
defined over $\oQ$.
\end{defi}

\begin{defi} \label{linear1}
  A linear subvariety $W \subset \nS_{\alpha}$ is a bi-algebraic
  subvariety $W= \pi(Y)$ such that the algebraic model of $Y$ is a (projectively)
  linear subspace $Z$ of $\PV_{\alpha}$. It is a $\oQ$-linear
  subvariety if moreover $W$ and $Y$ are
defined over $\oQ$.
\end{defi}

\begin{rem}
  Of course \Cref{bi-algebraic1} and \Cref{linear1}, given in terms of the
  developing map $D$, coincides with \Cref{bi-algebraic} and \Cref{linear}
  of the introduction, given in terms of local period charts.
\end{rem}

\subsection{Variants}
The content
of \Cref{subsec_strata} and \Cref{linear structure} holds for $\nH_\alpha$, with
the difference that the developing map will have values in
$\nV_\alpha$, thus defining an {\em integral linear structure} on
$\nH_\alpha^\an$.

As $\cS_\alpha$ is a finite \'etale cover of $\nS_\alpha$, the integral
projective structure on $\nS_\alpha$ induces one on $\cS_\alpha$, thus
defining a bi-algebraic structure on $\cS_\alpha$ making the finite
morphism $\cS_\alpha \to \nS_\alpha$ a morphism of bi-algebraic
structures in the obvious sense.

\subsection{Abelian varieties with many endomorphisms} \label{CM}
For the convenience of the reader we recall in this section classical results on (complex) abelian
varieties with many endomorphisms.

Let $A\simeq  A_1^{n_{1}}\times \cdots
\times A_k^{n_{k}}$ be a complex abelian
variety (where $\simeq$ denotes an isogeny and the $A_i$s are simple,
pairwise non isogenous). Let $\End(A)$ be its ring of endomorphism, and
$\End_\QQ(A):= \End(A) \otimes_\ZZ \QQ$. The $\QQ$-algebra
$\End_\QQ(A)$ is semi-simple, isomorphic to $\nM_{n_{1}}(D_1) \times
\cdots \times \nM_{n_{k}}(D_k)$, where $D_i:= \End_\QQ(A_i)$ is a
central division algebra over a finite extension $K_i$ of $\QQ$. We define the reduced degree of $\End_\QQ(A)$ by
\begin{equation} \label{reduced degree}
  [\End_\QQ(A): \QQ]_\red := \sum_{i=1}^k n_i [D_i:
  K_i]^{\frac{1}{2}}[K_i:\QQ]\;\; ,
\end{equation}
it is also the degree over $\QQ$ of a maximal \'etale subalgebra of $\End_\QQ(A)$.
As $\End_\QQ(A)$ acts faithfully on $H^1(A^\an, \CC)$, the reduced degree satisfies trivially:
\begin{equation} \label{inequality}
  2 \dim_\CC A \geq [\End_\QQ(A): \QQ]_\red\;\;,
\end{equation}

\begin{defi} 
The abelian variety $A$ is said to have complex multiplication (CM) if the
inequality~(\ref{inequality}) is an equality, or equivalently if
$\End_\QQ(A)$ contains an \'etale subalgebra of dimension $2
\dim_\CC A$, see \cite[Prop. 3.3]{Milne}.
\end{defi}

It follows immediately that
$A$ has CM if and only if each $A_i$ has CM, if and only if $D_i=K_i$ has degree $[K_i:\QQ] = 2 \dim_\CC
A_i$. With a little more effort, one shows that $K_i$ is a CM field
(i.e. a totally complex quadratic extension of a totally real field);
moreover if $A$ is isotypic then $A$ has CM if and only if $\End_\QQ(A)$ contains
a CM field of degree $2 \dim_\CC A$, see \cite[Prop. 3.6]{Milne}.

\begin{defi}
  The abelian variety $A$ is said to have real multiplication if
  $\End_\QQ(A)$ contains a totally real field $K_0$ of degree $[K_0:
  \QQ]= \dim_\CC A$.
\end{defi}

In particular, any abelian variety with complex multiplication has
real multiplication. One has the following partial converse:

\begin{lem}\label{lem_abvar}
Let $A$ be a complex abelian variety with real multiplication by $K_0$. Then
$A$ is isotypic, and $\End_\QQ A$ has a maximal \'etale subalgebra
which is either $K_0$, or a quadratic CM extension $K$ of $K_0$ (in
which case $A$ has CM).
\end{lem}

\begin{proof}
Let $g=\dim A$. Assume that
$A$ is not isotypic. Then there exists an isogeny $A \simeq A_1\times A_2$,
where the $A_i$s, $i=1,2$, are positive dimensional
abelian varieties without common simple factor up to isogeny. This realizes $K_0$ as a
subalgebra of $\End^0(A_i)$, $i=1,2$. The
numerical condition~(\ref{inequality}) forces $A_i$, $i=1, 2$, to have
CM with $\End_{\QQ} A_i= K_0$. As $K_0$ has a real place, it follows
from \cite[Lemma 3.7]{Milne} that $\dim_\CC A= [K_0:\QQ]$ divides
$\dim_\CC A_1$. This is a contradiction, thus $A$ is isotypic.

The field $K_0$ of degree $\dim_\CC A$ is contained in a maximal
\'etale subalgebra $E$ of
$\End_\QQ A$, which is of degree at most $2 \dim_\CC A$. Thus either $E=K_0$, or
$[E:\QQ]= 2\dim_\CC A$, in which case $A$ has CM, and $E$ can be
chosen to be a CM quadratic extension of $K_0$.
\end{proof}

\section{Characterization of arithmetic points} \label{SP:geometry}

  For $\alpha  =0$ (i.e. $C$ is an elliptic curve) the stratum
  $\nS_{0}$ coincides with the modular curve. \Cref{arithmetic
    points:geometry} says in this case that its arithmetic points are the CM points
  of the modular curve, a classical result of
  Schneider~\cite{Schneider37}. Schneider's theorem was generalized by
  Shiga and Wolfart \cite{SW} to say that if $A$ is a simple abelian variety
  over $\oQ$ admitting a non-zero algebraic one-form (defined
  over $\oQ$) all of whose periods are algebraic multiples one of each
other, then $A$ has complex multiplication. \Cref{arithmetic points:geometry} can
thus be thought as a generalisation ``to the mixed case'' of
\cite{SW}. As for most results nowadays in transcendence theory, our
main tool is W\"ustholz' Analytic Subgroup Theorem \cite[Theor.1 and
Corollary]{wus}:

\begin{theor} (W\"ustholz) \label{AST}
Let $\G$ be a connected commutative algebraic group over $\oQ$ with
Lie algebra $\fg$. Let $\exp: \fg \otimes_{\oQ} \CC \to \G(\CC)$ be the
exponential map for $\G(\CC)$.

Let $\fb \subset \fg $ be a positive dimensional $\oQ$-vector
subspace, $\BB \coloneqq \exp (\fb \otimes_{\oQ}\CC)$, and 
$u \in \fb \otimes_{\oQ}\CC$ such that $\exp(u) \in \G(\oQ)$. 
Let $\fh \subset \fb$ be the smallest $\oQ$-vector subspace such that
$u \in \fh \otimes_{\oQ} \CC$. Then $\fh$ is the Lie algebra of an algebraic
subgroup $\mathbf{H}$ of $\G$ defined over $\oQ$.
\end{theor}

\begin{proof}[\unskip\nopunct]
Let us now prove \Cref{arithmetic points:geometry}
Without loss of generality we can replace $\nS_\alpha$ by its finite
\'etale cover $\cS_\alpha$. Let us show that any arithmetic point $(C,[\omega]) \in
\cS_{\alpha}(\CC)$ satisfies necessarily the conditions stated in
\Cref{arithmetic points:geometry}. Condition~$(1)$ is obviously
satisfied. Let us choose $\omega \in H^0(C, \Omega^1_{C/ \oQ})$
representing $[\omega]$.

We first analyse the absolute periods of $\omega$ in order to prove
$(2)$. We denote by $\Alb(C)$ the Albanese of $C$ (a principally polarized abelian variety
defined over $\oQ$, isomorphic to the Jacobian of $C$). Poincar\'e's reduction theorem states that $\Alb(C)$ is
isogenous over $\oQ$ to a product $A_1^{k_{1}}\times \ldots \times
A_{N}^{k_{N}}$, with $A_\nu/\oQ$, $1\leq \nu \leq N$, pairwise
non-isogenous simple abelian varieties. Let 
$$\alb : C \to \Alb(C) \simeq \prod_{\nu =1}^N A_\nu^{k_{\nu}}$$ be the Albanese
map. It is uniquely defined up to fixing the image of a single point
in $C(\oQ)$. We normalize it by fixing $\alb(x_0)$ to be the identity
of $\Alb(C)$, where $x_0$ is the first zero of $Z(\omega)$ (this is where we use that we work
with the finite \'etale cover $\cS_\alpha$ rather than with $\nS_\alpha$). Let
$$ \alb_\nu: C \to A_\nu^{k_{\nu}}, \quad 1\leq \nu \leq N, $$ be the projections of
$\alb$ on each isogeny factor $A_\nu^{k_{\nu}}$. 

As the De Rham pull-back map 
$$ \prod_\nu \alb_\nu^*:  \prod_{\nu =1}^N {H^0(A^{k_{\nu}}_\nu , \Omega^1_{A^{k_{\nu}}_{\nu}/\oQ})} \to
  H^0(C, \Omega^1_{C/\oQ}) $$
is an isomorphism, there exists a
unique $\omega_\nu \in H^0(A^{k_{\nu}}_\nu, \Omega^1_{A^{k_{\nu}}_{\nu}/\oQ})$, $1 \leq
\nu \leq N$, such that
$\omega = \sum_{n=1}^N \alb_\nu^* \, \omega_\nu$. 

On the other hand the Betti pull-back map
$$\prod_\nu \alb_\nu^*:  \prod_{\nu} H^1((A_\nu^{k_{\nu}})^\an, \QQ) \to H^1(C^\an, \QQ)$$
is also an isomorphism and the comparison diagram between De Rham and
Betti cohomologies 
$$
\xymatrix{
\prod_{\nu} \left ( H^0(A^{k_{\nu}}_\nu, \Omega^1_{A^{k_{\nu}}_{\nu}/\oQ})  \otimes_{\oQ}
\CC \right) \ar[r]^<<<<{\prod_\nu \alb_\nu^*} \ar[d] & H^0(C,
\Omega^1_{C/\oQ}) \otimes_{\oQ} \CC
\ar[d] \\
 \prod_{\nu} \left( H^1((A^{k_{\nu}}_\nu)^\an , \QQ) \otimes_\QQ \CC \right ) \ar[r]_>>>>>{\prod_\nu \alb_\nu^*} \ar[r]
 & H^1(C^\an, \QQ) \otimes_\QQ \CC}
$$
commutes. Hence if the projectivized Betti class
$[\omega]_\Betti \in \proj H^1(C^\an, \CC)$ lies in $\proj H^1(C^\an,
\oQ)$, all the classes $[\omega_\nu]_\Betti \in \proj H^1((A^{k_{\nu}}_\nu)^\an,
\CC)$ corresponding to the non-zero factors $\omega_\nu$ lie in $\proj H^1((A^{k_{\nu}}_\nu)^\an,
\oQ)$.

It follows from \cite[Prop.3]{SW}  that for each $\nu$, $1 \leq \nu \leq N$, such that
$\omega_\nu \not =0$, the factor $A_\nu$ has complex
multiplication by a CM field $K_\nu$.

Moreover, it follows from \cite[prop.1 (1)]{SW} that there exists at
most one index $\nu$ such that $\omega_{\nu}\not =0$. With the
notation of \Cref{degree}, 
$A_{[\omega]}:= A^{k_{\nu}}_\nu$ has complex multiplication, thus,
in view of \Cref{CM}, satisfies the condition~$(2)(a)$ in
\Cref{arithmetic points:geometry}.

\medskip
Let us show Condition~$(2)(b)$ of \Cref{arithmetic
  points:geometry}.
Suppose first that $A_{[\omega]}$ is simple.
Let $K$ be the CM field of degree $2 d_{[\omega]}$ acting
on $A_{[\omega]}$. The homology $H_{1}(A_{[\omega]}^\an, \ZZ)$ is a
degree one module $\cO_K\cdot \gamma_0$ under an order $\cO_K$ of $K$. Let
us choose a basis of $H^0(A_{[\omega]}, \Omega^1_{A_{[\omega]}})$ over
  $\oQ$ consisting of eigendifferentials
$\omega_1, \dots, \omega_{d_{[\omega]}}$ for the action of $K$, with
  CM-type $\Phi = (\Phi_1, \cdots \Phi_{d_{[\omega]}})$. For each $1\leq i \leq d_{[\omega]}$, the
    periods $\int_\gamma \omega_i$, $\gamma \in
    H_{1}(A_{[\omega]}^\an, \ZZ)$, generate a one dimensional $\oQ$-vector
    space $\oQ\cdot \int_{\gamma_0} \omega_i$. Moreover, these
    $\oQ$-lines $\oQ\cdot \int_{\gamma_0} \omega_j$, $1\leq j \leq d_{[\omega]}$,
    are $\oQ$-linearly independent, see \cite[Prop. 1(2)]{SW}. Let us decompose
    $\omega_{[\omega]}= \sum_{i \in I \subset \{1, \cdots, d_{[\omega]}\}} \alpha_i \omega_i$, with
    $0 \not =\alpha_i \in \oQ$. The condition that the periods $\int_\gamma \omega$, $\gamma \in
    H_{1}(A_{[\omega]}^\an, \ZZ)$, generate a one dimensional $\oQ$-vector
    space is equivalent to saying that for all $k \in \cO_K$, the
    number $\Phi_i(k)$ is independent of $i \in I$. This implies that
    $|I|=1$, that is, $\omega_{[\omega]}$ is a
    $K$-eigendifferential. This finishes the proof of $(2)(b)$ in case
    $A_{[\omega]}$ is simple. In the general isotypic case
    $A_{[\omega]}= B^m$, choosing a $K$-eigenbasis for $H^0(B,
    \Omega^1_B)$ and taking for basis of
    $H^0(A_{[\omega]},\Omega^1_{A_{[\omega]}})$ the one obtained by
    pullback of $H^0(B,
    \Omega^1_B)$ under the $m$ projections, one easily checks that the
    periods of $\omega \in H^0(A_{[\omega]},\Omega^1_{A_{[\omega]}})$ generate a
    $\oQ$-vector space of dimension $1$ if and only if $\omega$ is a
    $K$-eigendifferential.

\sspace
Let us now analyse the relative periods of $\omega$ for proving
Condition $(3)$ of \Cref{arithmetic
  points:geometry},
namely that $\alb_{[\omega]}(Z(\omega))$ is a set of torsion 
points of $A_{[\omega]}$. Let us assume for simplicity that
$A_{[\omega]}$ is simple (the proof immediately adapts to the general
case) and let us just write $A$ for
$A_{[\omega]}$, $\omega_A$ for $\omega_{[\omega]}$ and $d$ for $d_{[\omega]}$. We denote by $\exp_{A}:
\Lie(A) 
\otimes_{\oQ} \CC \to A^\an$ the exponential for $A^\an$.

Choose $\omega_1:= \omega_{A}, \omega_2, \ldots, \omega_d$ a basis of $H^0(A,
\Omega^1_{A/\oQ})$ completing $\omega_A$. Let $\omega_1^*, \ldots,
\omega_d^*$ be the corresponding dual basis of $\Lie(A)$. Let $$p: H_1(A^\an;
\CC) \to H^{-1, 0}(A) = \Lie(A) \otimes_{\oQ} \CC$$ be the natural
Hodge projection, it associates to any element $\gamma \in H_1(A^\an; \CC)$ the vector
$p(\gamma)$ with coordinates 
$(\int_\gamma \omega_1, \ldots ,\int_\gamma \omega_d)$ in the basis
$(\omega_i^*)$. The map $p$ defines an isomorphism of $H_1(A^\an,
\ZZ)$ with the period lattice $\Gamma:= p(H_1(A^\an,
\ZZ) ) \simeq \ZZ^{2d} \subset \Lie(A) \otimes_{\oQ} \CC$. Moreover, as
$A$ has complex multiplication by a number field $K$, there exists a unique $\oQ$-vector
subspace $V_{\oQ} \subset \Lie(A) \otimes_{\oQ} \CC$ of dimension
$d$ containing $\Gamma$.

\begin{rem}
In other words: $V_{\oQ}$ defines a new $\oQ$-structure on $\Lie(A)
\otimes_{\oQ} \CC$.
\end{rem}

Let $x$ be a point of $Z(\omega)$, different from $x_0$. Let $\beta $ be a path in $A$
joining $x_0$ to $x$, and $\ti{x} \in \Lie(A)
\otimes_{\oQ} \CC$ the corresponding element. Hence $\exp(\ti{x})=
\alb_{[\omega]}(x)$ and $\ti{x} = (\int_\beta \omega_A, \int_\beta \omega_2, \ldots,
\int_\beta
\omega_d)$. 

\begin{lem}
$\ti{x} \in V_{\oQ}$.
\end{lem}

\begin{proof}
We follow the strategy of the proof of \cite[Prop.3]{SW}.

Let us fix a non-zero $\gamma \in H_1(A^\an, \ZZ)$. As
noticed in \cite[proof of Proposition 1]{SW}, the simplicity of $A$
and \Cref{AST} imply that $\int_\gamma \omega \not =0$.
The condition $[\omega]_{\Betti} \in \proj H^1(C,
Z([\omega]); \oQ)$ implies that there exists $\alpha \in \oQ$ such
that $\int_\beta \omega_A = \alpha \cdot \int_\gamma \omega_A$.

Consider $\exp_{A\times A} : \Lie(A \times A)\otimes_{\oQ} \CC \to (A
\times A)^\an$, the exponential for $(A \times A)^\an$. Let $\fb
\subset \Lie(A \times A)$ the $\oQ$-hyperplane defined by 
$$ \fb:= \{(z, w) \in \Lie(A) \times \Lie(A) \;\;| \;\; \omega_A(w) =
\alpha \cdot \omega_A(z) \;\}\;\;.$$
The vector $u:=(p(\gamma), \ti{x}) \in \fb \otimes_{\oQ}\CC$ is non-zero and satisfies $\exp_{A\times A}(u)=
(1, x) \in (A \times A)(\oQ)$. Applying \Cref{AST} to $\G:= A \times
A$, we conclude that there exists a non-zero $\oQ$-subspace $\fh
\subset \fb$ which is the Lie algebra of an algebraic
subgroup $\mathbf{H}$ of $A \times A$ defined over $\oQ$. As $A$ is
simple, the group $\mathbf{H}$ is isogenous to $A$. By \cite[Lemma 3]{SW}
the Lie algebra $\fh$ is defined by a non-singular system of $d$
linear equations $w = C \cdot z$ for some $C \in \End_0(A)$. As $A$ is
simple, with complex multiplication by a number field $K$, $\End_0(A)=
K$. Hence $C = \alpha \in K$ and $\ti{x} = \alpha \cdot p(\gamma) \in V_{\oQ}$.
\end{proof}

We conclude that $x$ is a torsion point in $A$ by applying the
following result of Masser \cite{Ma76} (a similar statement was 
obtained by Lang \cite{Lang75}). 

\begin{theor} (Masser) \label{Masser}
Let $A$ be an abelian variety of dimension $r$ with CM (in particular
$A$ is defined over $\oQ$). Let $\exp_A : \Lie(A_\CC) \to A^\an$ be
the uniformization for $A^\an$, with kernel $\Gamma$. Let $V_{\oQ}
\subset \Lie(A_\CC)$ be the $\oQ$-vector subspace of dimension $r$
containing $\Gamma$.

Let $ \ti{x} \in V_{\oQ}$ such that $x:= \exp_A (\ti{x})$ belongs to
$A(\oQ)$. Then $x$ is a torsion point in $A$ and $\ti{x} \in \Gamma
\otimes_\ZZ \QQ$.

\end{theor}

This finishes the proof that any arithmetic point $(C, [\omega]) \in
\nS_{\alpha}(\CC)$ satisfies $(a)$, $(b)$ and $(c)$ of
\Cref{arithmetic points:geometry}.

The same analysis shows easily that the converse is true. This
finishes the proof of \Cref{arithmetic points:geometry}.
\end{proof}

\begin{rem}
Masser's \Cref{Masser} is itself an incarnation of the bi-algebraic
  format over $\oQ$ for abelian varieties with complex multiplication,
  the algebraic structure over $\oQ$ considered on $\Lie(A_\CC)$ being
  the one given 
by the period lattice $\Gamma$ rather 
than the standard one $\Lie A$. See \cite[Section 2.2.2]{Ull16}.
It is also a direct consequence of \Cref{AST}. See \cite[Theor. 2.9]{Ull16} for a
generalization.
\end{rem}


\section{Distribution of arithmetic points}

\subsection{Density of arithmetic points: proof of \Cref{arithmetic
    points:density}} \label{density}

\begin{proof}[\unskip\nopunct]
For $\alpha=0$, the stratum $\nS_0$ is the modular curve parametrizing
complex elliptic curves. The arithmetic points of $\nS_0$ are
precisely the CM elliptic curves, which are
well known to be dense in $\nS_0^\an$, and are obviously of degree~$1$.

For $\alpha \not = 0$, we construct a
 dense set of arithmetic points of degree~$1$ in $\nS_\alpha^\an$ as covers of CM
 elliptic curves ramified precisely over torsion points. Let us notice
 we do not need \Cref{arithmetic points:geometry} for this construction.

Let $U\subset \nS^\an_\alpha$ be a simply-connected open subset and
$\phi: U\to \proj^{d_\alpha}(\CC)$ be a period
chart. Let us show that all points in $U$ with homogeneous $\phi$-coordinates in $\QQ(i)
\subset \CC$ are arithmetic of degree~$1$
(they are clearly dense in $\nS_\alpha^\an$). Let $u\in U$ be such a
point, corresponding to a pair $(C,[\omega])$. Choose a representative
$\omega$ of $[\omega]$ such that its period coordinates lie in
$\QQ(i)$. Let $\gamma_0,\ldots, \gamma_{d_\alpha}$ be a basis of
$H_1(C^\an,Z(\omega);\ZZ)$ and let $D\in\NN$ be an integer such that
$\int_{\gamma_k} \omega\in \frac{1}{D}(\ZZ+i\ZZ)$, for $0\le k\le
{d_\alpha}$. If $x_0\in Z(\omega)$ and $E$ is the complex elliptic
curve defined by $E^\an=\CC/\frac{1}{D}(\ZZ+i\ZZ)$, the map 
\begin{align*}
f:C^\an&\to E^\an\\
x&\mapsto \int_{x_0}^x \omega
\end{align*}
defines a ramified cover, whose ramification locus lies over $0\in
E^\an$. Moreover $f^*(dz)=\omega$,  where $dz$ is the 1-form on $E$ 
deduced from $dz$ on $\CC$. The elliptic curve $E$ is defined over
$\QQ$ (indeed its analytification is isomorphic to $\CC/(\ZZ+i\ZZ)$,
which has $j$-invariant $1728$, hence defined over $\QQ$), in
particular over $\Qbar$. Since $0\in E(\Qbar)$ and $C$ is ramified
only over $0$, we deduce by \cite[XIII, Proposition 4.6]{SGA1} that
$C$ and $f$ are defined over $\Qbar$ (in the notation of {\em
  loc. cit.} one takes $k=\Qbar$ and $Y=\Spec(\CC)$). Moreover, the
line $\CC\omega$ in $H^0(C,\Omega^1_C)$ is defined over $\Qbar$: indeed,
since $H^0(E,\Omega^1)$ is one dimensional, there exists
$\alpha\in\CC$ such that $\alpha \,dz\in H^0(E_\Qbar,\Omega^1_{E/\oQ})$, and
then $\alpha \, \omega$ is defined over $\Qbar$. 
\end{proof}

\begin{rem}
(1) An abelian differential $(C,\omega)$ with relative periods in
$\QQ(i)$ is known in the literature as a {\em square-tiled surface}. 
  
(2) The above argument goes through with $\QQ(i)$ replaced by any
quadratic imaginary number field embedded in $\CC$.
\end{rem}

Let us describe an apparently more general construction of arithmetic
points of degree $1$, by allowing
several branch torsion points in $E$ rather than just one. Let $E$ be a complex elliptic curve with complex multiplication (we
refer to \cite[Chap.II]{Sil} for a description of such curves). Hence $E$
admits a model over $\oQ$, which we still denote $E$ by abuse of
notations. Let us write its complex uniformization $E^\an= \CC /
(\ZZ + \tau \ZZ)$, where $\tau$ is imaginary quadratic. Let 
$\omega_E \in
H^0(E, \omega^1_{E/\oQ})$
be a non-zero algebraic one-form on $E$
defined over $\oQ$. Its pull-back to $\CC$ under the uniformization
map $\pi: \CC \to E^\an$ can be written
$\tilde{\omega} = \alpha \cdot dz$, where $z$ is a natural holomorphic
coordinate on $\CC$. The absolute periods of $\omega_E$ on
$H_1(E^\an;\ZZ) = \ZZ + \tau \ZZ$ lie in $\alpha\cdot (\ZZ + \tau 
\ZZ)$ hence in $\alpha\cdot \oQ$.

Let $k$ be an integer and let $e_i \in E(\oQ)$, $1\leq i
\leq k$, be distinct torsion points in $E$. Choose $\tau_i\in \CC$, $1 \leq
i\leq k$, some $\pi$-preimage of $e_i$ under the uniformization
map $\CC \to E^\an$. If $n_i$ denotes the order of $e_i$ in $E$, the complex
number $n_i \cdot \tau_i $ lies in the lattice $\ZZ + \tau \ZZ$. In
particular 
\begin{equation} \label{alpha}
\forall i\in\{1, \ldots, k\}, \qquad \int_0^{\tau_i} \tilde{\omega}_E\in \alpha \cdot
\oQ\;\;.
\end{equation}
This shows that
\begin{equation} \label{omegaE}
[\omega_E]_\Betti \in \proj H^1(E^\an, \{e_1, \ldots, e_k\};
\oQ)\;\;.
\end{equation}

Let $g\geq 2$ and choose a type $\alpha$ for $g$. Given $f:C \to E$ a smooth
projective curve over $\oQ$ ramified only 
over the $e_i$'s, we define $\omega:= f^*\omega_E \in H^0(C, \Omega^1_{C/\oQ})$ as
the pull-back to $C$ of $\omega_E$. Classical ramification theory
shows one can choose $k$ and $f$ such that the pair $(C, [\omega])$ defines a point in
$\nS_{\alpha}^\an$. As the point $[\omega] \in \proj H^1(C^\an,
Z(\omega) ; \CC)$ is the image of $[\omega_E]$ under the complexification of
the natural Betti pull-back map $$f^*:  \proj H^1(E^\an, \{e_1, \ldots, e_k\}; \ZZ) \to  \proj
H^1(C^\an, Z(\omega) ; \ZZ) \;\;,$$ it follows from (\ref{omegaE})
that $[\omega]_\Betti$ belongs to $\proj H^1(C, Z(\omega); \oQ)$.
Hence the point $(C, [\omega])$ is arithmetic in $\nS_{\alpha}(\CC)$.

\medskip
Let us notice that the construction of arithmetic points of degree $1$
we just described does not produce more points in $\nS_\alpha$ than
square-tiled surfaces. Indeed, if $f: C\to E$ is a ramified cover
whose whose branch locus consists of torsion points, let $N \in \NN$
be an exponent killing these torsion points and $[N]: E \to E$ be the
multiplication-by-$N$ map. Then $[N] \circ f: C\to E$ is a ramified
cover, ramified only over $0 \in E$. The abelian differentials $(C,
f^* \omega_E)$ and $(C, ([N] \circ f)^* \omega_E)$ are different
points of $\nH_\alpha$, but define the same point in $\nS_\alpha$.

\subsection{Propagation of arithmetic points by ramified
  covers} 

The construction~\Cref{density} can be generalized to the following
lemma, whose proof using \Cref{arithmetic
  points:geometry} is immediate:

\begin{lem} \label{subsec_constr_points}
Let $(C,[\omega])$ be an arithmetic point in $\nS_\alpha(\oQ)$, and
let $\alb_{A_{[\omega]}}:C\to {A_{[\omega]}}$  be as in \Cref{degree},
normalised by sending a zero of
$\omega$ to the origin in ${A_{[\omega]}}$. If $f:C'\to C$ is a ramified cover of curves whose branch locus
$\textrm{Branch}(f)$ in $C$ satisfies $\alb_{A_{[\omega]}}(\textrm{Branch}(f))\subset
\textrm{Tors}({A_{[\omega]}})$, then $(C',[f^*\omega])$ is also an arithmetic point
(in the corresponding stratum $\nS_{\alpha'}(\oQ)$).

Moreover, the degrees of $(C,[\omega])$ and $(C',[f^*\omega])$ are equal.
\end{lem}

\begin{proof}
  With our normalisation, points of $Z(\omega)$ are sent to torsion
points in ${A_{[\omega]}}$, by condition \Cref{arithmetic points:geometry}(c).  
The variety ${A_{[f^*\omega]}}$ factor of $\Alb(C')$ coincides with
$A_{[\omega]}$, the image of $C'$ in ${A_{[\omega]}}$ is $\alb_{A_{[\omega]}}(C)$ and the 
zeros of $f^*\omega$ lie over torsion points. Hence $(C',[f^*\omega])$
is an arithmetic point in $\nS_{\alpha'}(\oQ)$ by \Cref{arithmetic
  points:geometry}, of the same degree as $(C,[\omega])$.
\end{proof}

\subsection{Primitive arithmetic points}

\begin{defi} \label{geomPrim}
  An arithmetic point $(C, [\omega]) \in \nS_{\alpha}(\oQ)$ is said to be geometrically primitive if it
  cannot be obtained from an arithmetic point in a stratum of smaller
  genus via a ramified cover as in \Cref{subsec_constr_points}.
\end{defi}

\begin{rem}
Notice that any geometrically primitive arithmetic point
$(C,[\omega])$ in $\nS_\alpha(\oQ)$ has degree $d_{[\omega]}$ at least $2$ as soon as
$\alpha \not=0$.
\end{rem}

\begin{defi} \label{algPrim}
 An arithmetic point $(C, [\omega]) \in \nS_{\alpha}(\oQ)$ is said to
 be algebraically primitive if its associated CM factor $A_{[\omega]}$ coincides
 with $\Alb(C)$, i.e. $d_{[\omega]}= g(C)$.
\end{defi}

\begin{rem}
 Any algebraically primitive arithmetic point $(C, [\omega]) \in
 \nS_{\alpha}(\oQ)$ is obviously geometrically primitive. 
\end{rem}

\begin{Example}
Let us provide examples of algebraically primitive arithmetic points of degree at least $2$
which are Veech surfaces. The simplest examples of algebraically primitive Teichm\"uller curves 
belong to the original Veech's family \cite[Theorem 1.1]{Ve89}. For all $n\ge 1$, these Teichm\"uller curves
are the projection in $\proj\Omega\cM_{\lfloor{\frac{n-1}{2}}\rfloor}$
of the $\GL^+(2,\RR)$-orbit of $(C_n, \frac{dx}{y})$
where $C_n$ is the
(smooth projective model of the) plane hyperelliptic curve
$y^2=x^n-1$. When $n$ is odd, they belong to the minimal 
stratum $\nS_{n-3}$ of $\proj\Omega\cM_{\frac{n-1}{2}}$ parametrizing
abelian
differentials with a single zero. When $n=p$ is prime, it is
well-known that $\Jac(C_p)$ is CM and simple; since
$\omega_p$ has a single zero, \Cref{arithmetic points:geometry}
implies that $(C_p,[\omega_p])$ is an algebraically primitive
arithmetic point in $\nS_{p-3}$ which is a Veech surface.

\medskip
Similarly, let $C_{1,1}$ be the genus $2$ hyperelliptic curve $y^2 =
x^6-x$. By \cite[Remark 13.3.8]{BL} its Jacobian surface $\Jac(C_{1,1})$ has
complex multiplication. Moreover it is simple: according to Bolza
\cite{Bolza} (see also \cite[p.340]{BL}) its reduced group of
automorphism is $\ZZ/5\ZZ$, hence $C$ does not admit any involution
with quotient an elliptic curve. The abelian differential $(C_{1,1},
[\frac{dx}{y}])$ is thus an algebraically primitive arithmetic point
in $\nS_{1,1}(\oQ)$ which is a Veech surface. The algebraically primitive Teichm\"uller curves it generates is the
famous decagon family, see \cite[Theor. 1.1]{McM06} (we thank
M. M\"oller for mentioning this example to us).
\end{Example}

\begin{rem} \label{bow}
Veech's examples are a particular case of a more general family of
Teichm\"uller curves discovered by Bouw and M\"oller \cite{BM10},
whose Veech groups are triangle groups with parameters $(m,n,\infty)$,
$m,n\in \NN\cup\{\infty\}$. The Bouw-M\"oller family also contains
algebraically primitive Teichm\"uller curves in unbounded genera,
different from Veech's. It is at present the only known
family containing primitive Teichm\"uller curves in unbounded genera.

We claim that every Teichm\"uller curve in the Bouw-M\"oller family contains an
arithmetic point. Indeed let $\mathscr{H}\to \overline{C}$ be the
family of curves constructed in \cite[Section 6]{BM10}. Thus $\overline{C}\xrightarrow{\pi}\proj^1$ is a suitable cover
ramified over $\{0,1,\infty\}$, such that the image of
$\overline{C}\setminus \pi^{-1}(\infty)$ in $\cM_g$ is a Teichm\"uller
curve. If $c\in\pi^{-1}(0)$, one shows that the fiber $\mathscr{H}_c$ has a
Jacobian with CM. This follows by inspection of the construction in \cite[Lemma 6.8]{BM10}
together with the fact that (a connected component of) the curve
$z^N=x^a(x-1)^b$ has Jacobian with CM.  By \Cref{arithmetic points:geometry} and the fact
that $\mathscr{H}_c$ together with its 1-form up to scaling $[\omega]$
is a Veech surface, we deduce that $(\mathscr{H}_c,[\omega])$ is an
arithmetic point. 
\end{rem}

\subsection{Proof of \Cref{arithm-versus-Veech}} \label{ex_prim_Teich}

\begin{proof}[\unskip\nopunct]

Notice that the statement of \Cref{subsec_constr_points} remains valid
if one replaces ``arithmetic point'' by ``Veech surface''. As any
arithmetic point $(C, [\omega]) \in \nS_\alpha$ of degree $1$ is
obtained as in \Cref{subsec_constr_points} from a CM elliptic curve
$(E, [\omega_E]) \in \nS_0$, and as $(E, [\omega_E])$ is trivially a
Veech surface in $\nS_0$, it follows that any arithmetic point of
degree $1$ is a Veech surface.

\medskip
Let us exhibit an arithmetic point of higher degree which is not a Veech surface.
The curve $C_7: y^2=x^7-1$ together with the differential
$\left[\frac{x\,dx}{y}\right]$ defines an arithmetic point of degree $3$. To see
this we check the conditions of \Cref{arithmetic points:geometry}. As
remarked above, $\Jac(C_7)$ is simple of dimension $3$ and has complex
multiplication. Denote by $\infty$ the unique point of $C_7$ that lies
above $\infty\in\proj^1$ for the projection $(x,y)\mapsto x$. The
divisor of zeros of  $\left[\frac{x\,dx}{y}\right]$ is given by
$(0,i)+(0,-i)+2\cdot\infty$. The classes $[(0,\pm i) - \infty]$ are
7-torsion in the Jacobian (as $7(0,i) - 7\infty$ is the divisor of the
function $y-i$), hence condition of (c) of \Cref{arithmetic
  points:geometry} is fulfilled. 
However, by \cite[Remark 2.9]{Moe06} and \cite[Theorem 7.5]{McM03},
this abelian differential does not generate a Teichm\"uller curve. 
\end{proof}


\subsection{Linear invariant subvarieties with only few arithmetic points: proof of
  \Cref{prop_counterexample}}\label{sec_counterexample} 

\begin{proof}[\unskip\nopunct]

Let us first give the proof of \Cref{prop_counterexample} in the case
of a Teichm\"uller
curve, where the geometry is simpler. Recall the following fundamental result:
\begin{theor}[\cite{Moe06}]\label{thm_moe06}
Let $Z\subset\proj\Omega^1\cM_g$ be a Teichm\"uller curve of
degree~$r$. Then the restriction of $\bV_\alpha$ to $Z$ satisfies 
$$\Gr_1^W {\bV_\alpha}_{|Z}= (\bigoplus_{i=1}^d \bL_i) \oplus \bM\;\;,$$
where $\bM$ is defined over $\QQ$, the $\bL_i$ are rank $2$ local subsystems, maximal Higgs for $i=1$,
and non-unitary but not maximal Higgs for $i \not =1$.

In particular, 
the image ${Z'}^\pure$ of $Z$ in $\cA_g$ is contained in the locus of
abelian varieties 
that split up to isogeny as $A_1 \times A_2$, $\dim A_1 = d$ and $A_1$ has real
multiplication. 
\end{theor}

\noindent Let $Z\subset \nS_\alpha \subset \proj\Omega^1\cM_g$ be a Teichm\"uller curve of degree
$d$. Consider the following diagram provided by (\ref{defi period
  map}) and \Cref{thm_moe06}:
\begin{equation} \label{diag1}
  \xymatrix{
  \nS_\alpha \ar@{^(->}[r]^\Phi & \fA_g^{(n-1)} \ar@{->>}[r] &\cA_g \\
  Z \ar@{^(->}[u] \ar@{^(->}[r] \ar@{->>}[d] & \fA_{d, g-d}^{(n-1)} \ar@{->>}[r] \ar@{->>}[d] \ar@{^(->}[u] &\cA_r
  \times \cA_{g-d} \ar@{->>}[d] \ar@{^(->}[u] \\
  Z' \ar@{^(->}[r] &\fA_{d}^{(n-1)} \ar@{->>}[r] &\cA_d.
  }
\end{equation}

\noindent
Here $\fA_{d, g-d}^{(n-1)} \subset \fA_g^{(n-1)}$ is the special
subvariety of $\fA_g^{(n-1)}$ containing $Z$ and parametrizing abelian varieties 
that split up to isogeny as $A_1 \times A_2$, $\dim A_1 = d$. It is
isomorphic to $\fA_{d, g-d} \times_{(\cA_d
  \times \cA_{g-d})} \cdots \times_{(\cA_d
  \times \cA_{g-d})} \fA_{d, g-d} $, and
$\fA_{d}^{(n-1)}: = \fA_d \times_{\cA_d} \cdots \times_{\cA_d} \fA_d$
is its natural quotient.

Let ${Z'}^\pure \subset \cA_d$ be the projection of $Z'$ in $\cA_d$
(its pure part).
If $(C,[\omega])\in Z(\CC)$ is an arithmetic point, it follows from \Cref{arithmetic points:geometry}
that the image of $(C,[\omega])$ in $Z''(\CC)\subset \cA_d(\CC)$ is a CM point. Therefore, if $Z$ contains
infinitely many arithmetic points, then $\overline{{Z'}^\pure}$ contains
infinitely CM points. By the Andr\'e-Oort conjecture for $\cA_d$
\cite{Ts18}, $\overline{{Z'}^\pure} \subset \cA_d$ is a
Shimura curve in $\cA_d$.

\noindent
It then follows from \cite[Theor.1.2]{MolST} that the $\bL_i$s, $2
\leq i \leq d$, which comes from the Shimura curve $\overline{{Z'}^\pure}$, are also
maximal Higgs. This is a contradiction to \Cref{thm_moe06} if $d\geq 2$ . This proves
\Cref{prop_counterexample} for Teichm\"uller curves.


\medskip
Let us now turn to the general case. Let $Z \subset \nS_\alpha$ be
a linear invariant subvariety. Recall the variation~(\ref{vmhs}) of mixed $\ZZ$-Hodge structures
($\ZZ$VMHS) on
$\nS_\alpha$: 
\begin{equation} 
  0 \to W_0 \bV_{\alpha, \ZZ} \simeq \ZZ(0)^{n-1} \to \bV_{\alpha, \ZZ} \to
  \Gr_1^W \bV_{\alpha, \ZZ}  \to 0\;\;.
  \end{equation}
The real tangent bundle 
$T_\RR Z$ is a sub-$\RR$VMHS of the
restriction to $Z$ of $\bV_{\alpha, \RR}:= \bV_{\alpha, \ZZ}
\otimes_\ZZ \RR$, of the form:
\begin{equation} \label{funda}
  0 \to W_0 T_\RR Z \to T_\RR Z \to \bH_1: = \Gr_1^W T_\RR Z \to
  0\;\;. 
  \end{equation}

  \noi
 We denote by $u$ the rank of $W_0 T_\RR Z$ and by $2k$ the rank of
 $\bH_1$. The integer $k$ is called the rank of $Z$. In particular:

  \begin{equation} \label{dim1}
    \dim_\CC Z:= u +2k -1\;\;.
    \end{equation}

According to \cite{Fi}, the $\RR$VMHS $T_\RR Z$ is defined over a
totally real number field $K_0 \stackrel{\iota_1}{\hookrightarrow} \RR$:
$$ T_\RR Z = T_{K_{0}} Z \otimes_{K_{0}, \iota_1} \RR\;\;.$$
We write $d:= [K_0:\QQ]$. In particular 
\begin{equation} \label{splitting}
 \Gr_1^W {\bV_{\alpha,\RR}}_{|Z}= (\bH_1 \oplus
 \cdots \oplus \bH_d) \oplus \bM\;\;,
\end{equation}
where the $\bH_i$, $1 \leq i \leq d$, are the
$\Gal(\oQ/\QQ)$-conjugates of $\bH_1$ and the factor $\bM$ is defined over $\QQ$. We thus obtain the following
generalization of (\ref{diag1}) (Teichm\"uller curves correspond to
the case $(k=1, u=0)$):

\begin{equation} \label{diag2}
  \xymatrix{
  \nS_\alpha \ar@{^(->}[r]^\Phi & \fA_g^{(n-1)} \ar@{->>}[r] &\cA_g \\
  Z \ar@{^(->}[u] \ar@{^(->}[r] \ar@{->>}[d] & \fA_{kd, g-kd}^{(n-1)}
  =\fA_{kd}^{(n-1)} \times \fA_{g-kd}^{(n-1)}  \ar@{->>}[r] \ar@{->>}[d] \ar@{^(->}[u] &\cA_{kd}
  \times \cA_{g-kd} \ar@{->>}[d] \ar@{^(->}[u] \\
  Z' \ar@{^(->}[r] &\fA_{kd}^{(n-1)} \ar@{->>}[r] &\cA_{kd}.
  }
\end{equation}

\begin{lem} \label{iso}
The projection $Z \to Z'$ is an
isomorphism. In particular $\dim_\CC Z = \dim_\CC Z'$.
\end{lem}

\begin{proof}
Let $(C, [\omega]) \in Z$. Seen in $\fA_g^{(n-1)}$, this point corresponds to
the mixed Hodge structure
\begin{equation} \label{eq1}
  0 \to \ZZ(0)^{n-1} \to H^1(C, Z([\omega]);\ZZ) \to H^1(\Alb(C), \ZZ) \to
  0\;\;.
  \end{equation}
Its image in $Z' \subset \fA_{kd}^{(n-1)}$ is the mixed Hodge structure
\begin{equation} \label{eq2}
  0 \to \ZZ(0)^{n-1} \to E \to H^1(A_{[\omega]}, \ZZ) \to
0 \end{equation} 
corresponding to the projection of the zeroes of $[\omega]$ in
$A_{[\omega]}$. But (\ref{eq1}) is nothing else than the pull-back of (\ref{eq2})
under the natural projection $\Alb(C) \to A_{[\omega]}$. Hence the result.
\end{proof}

Suppose that $k=1$ and $(C,[\omega])\in Z(\CC)$ is any arithmetic
point, or that $k>1$ and $(C,[\omega])\in Z(\CC)$ is an arithmetic point of degree at least
$kd$. In all these cases the containment $A_{[\omega]}\in \cA_{kd}$
implies that $(C,[\omega])$ has degree exactly $kd$. It then follows from \Cref{arithmetic points:geometry}
that the image of $(C,[\omega])$ in $Z'(\CC)\subset
\fA_{kd}^{(n-1)}(\CC)$ is a CM point. Therefore, if $Z$ contains 
a Zariski-dense set of arithmetic points of such points, then
$\overline{Z'}^\Zar \subset \fA_{kd}^{(n-1)}$ contains 
a Zariski-dense set of CM points. By the Andr\'e-Oort conjecture for the mixed
Shimura variety of Kuga type $\fA_{kd}^{(n-1)}$
\cite{Ts18}, \cite{Gao}, the closed subvariety $\overline{Z'}^\Zar
\subset \fA_{kd}^{(n-1)}$, which contains $Z'$ as a Zariski-open dense
subset, is a special subvariety of $\fA_{kd}^{(n-1)}$.

In \cite[Cor. 1.7]{Fil17}, extended in
\cite[Prop. 4.7]{EFW18}, it is shown that the real
Zariski-closure of the monodromy of the flat bundle $T_\RR Z$ at a
point $z_0 \in Z$ is as
big as it can be, namely equal to $\Sp(\bH_{1, z_{0}}) \ltimes \HHom
(\bH_{1, z_{0}}, (W_0 T_\RR Z)_{z_{0}})$; and that the same holds true for the
Galois conjugates bundles: the real Zariski closure of $(T_\RR Z)_{\iota_{i}}$,
$1 \leq i \leq d$, is $\Sp(\bH_{i, z_{0}}) \ltimes \HHom
(\bH_{i, z_{0}}, (W_0 T_\RR Z)_{z_{0}})$. If we denote by $T_\QQ Z$
the $\QQ$-sub-VMHS $\Res_{K_0/\QQ} T_{K_{0}} Z$ of $\bV_{\alpha, \QQ}$,
it follows that the ($\QQ$-)algebraic monodromy group of $T_\QQ Z$ is
the group $\Res_{K_0/\QQ} (\Sp_{2k, K_0} \ltimes \HHom_{K_{0}}(\bH_{1,
  K_{0},  z_{0}}, (W_0 T_\RR Z)_{K_{0}, z_{0}}))$. In particular the 
$\QQ$-algebraic monodromy group defining the special subvariety $\overline{Z'}^\Zar$
of $\fA_{kd}^{(n-1)}$ equals $\Res_{K_0/\QQ} (\Sp_{2k, K_{0}} \ltimes \HHom_{K_{0}}
(\bH_{1, K_{0},  z_{0}}, (W_0 T_\RR Z)_{K_{0}, z_{0}}))$. It
follows that:
\begin{equation} \label{dim2}
  \dim \ol{Z'}^\Zar = d \times (\frac{k(k+1)}{2} + uk)\;\;.
\end{equation}

Comparing (\ref{dim1}) and (\ref{dim2}) using \Cref{iso}, one obtains:
\begin{equation} \label{ineq}
  u+2k-1 =d \times (\frac{k(k+1)}{2} + uk) \;\;.
\end{equation}

One easily checks that this inequality can be satisfied only for
$(d=1, k=2, u=0)$ or $(d=1, k=1, \textnormal{any}\, u)$.
This finishes the proof of \Cref{prop_counterexample}.

\end{proof}



\section{Bi-algebraic subvarieties of $\nS_{\alpha}$ and
  linearity} \label{linearity}

Let us introduce a few notations we will use in this section. Working
over $\nH_\alpha$ or $\nS_\alpha$, we write:

\begin{itemize}
\item $\bV^\pure_{\alpha,\ZZ}$, resp. $\bV^\pure_\alpha$, for
the pure piece of weight one $\Gr^W_1\bV_{\alpha,\ZZ}$, resp. $\Gr^W_1\bV_\alpha$
(either on $\nH_\alpha$ or $\nS_\alpha$, depending on the context); $\cV^\pure_\alpha$ for the holomorphic vector
bundle defined by $\bV^\pure_\alpha$; and $\nV^\pure_\alpha$ for the fiber of
$\bV^\pure_\alpha$ over a point.

\item $\cC_g\to \cM_g$ for the universal curve. 
\item $\cC_\alpha:=\pi_\alpha^*\cC_g$ and
  $\cJ_\alpha:=\Jac(\cC_\alpha/\nH_\alpha)$ the relative Jacobian of
  $\cC_\alpha$ over $\nH_\alpha$ or $\nS_\alpha$ (where $\pi_\alpha$
  denotes the projection from $\nH_\alpha$ or $\nS_\alpha$ to $\cM_g$).
\end{itemize}

Let $Z\subset\nS_\alpha$ be an algebraic subvariety. Let $\mathbb{T}$
be the largest constant local sub-system of $\Gr_1^W \bV_\alpha|_Z$,
$\cT:=\mathbb{T}\otimes \cO_Z$ and $\omega^\pure|_Z$ the projection of
$\omega|_Z$ to $\proj\cV^\pure_\alpha|_Z.$ Our analysis of the
bi-algebraic structure on $\nS_\alpha$ depends on the position of
$\omega^\pure|_Z$ with respect to $\proj\cT$. We deal with two
``opposite'' cases: 
\begin{itemize}
\item $\omega^\pure|_Z$ is a constant section of $\proj\cT$. 
\item $\omega^\pure|_Z$ is not generically contained in $\proj\cT$.
\end{itemize}

We do not deal with the intermediate case where $\omega^\pure|_Z$ is a non-constant section of $\proj\cT$.
The first case corresponds to the case where $Z$ is contained in an isoperiodic leaf. We treat it in \Cref{subsec_isoper}.
The second case is treated in \Cref{subsec_cond_star}. In
\Cref{subsec_gen_2} we prove that the intermediate case mentioned
above does not occur in genus 2 and in \Cref{genus3} that it
does occur in genus $3$.

\subsection{The isoperiodic foliation and proof of \Cref{bi-algebraic=linear1}}\label{subsec_isoper}

The isoperiodic foliation is a foliation on $\nH_\alpha$ by complex
submanifolds of codimension $2g$ defined locally as follows. Over a
chart $U\subset\nH_\alpha^\an$ with period coordinates $U\to
\nV_\alpha$, the isoperiodic foliation on $U$ is given as pullback of
the foliation on $\nV_\alpha$ given by the fibers of the projection
$\nV_\alpha\to \nV_\alpha^\pure$. In other words, a leaf around a point
$(C_0,\omega_0)$ consists of nearby points $(C,\omega)$ such that the
absolute periods of $\omega$ are the same as the absolute periods of
$\omega_0$ (the periods are computed on a flat frame of homology). By
passing to the projectivisations, one obtains the isoperiodic
foliation on $\nS_\alpha$.  


\begin{lem}\label{lem_alg_in_isoper}
Let $Z\subset\nH_\alpha$ be an algebraic subvariety of an isoperiodic leaf.  There exists an isogeny of abelian schemes
$$f: \cJ_{\alpha}|_Z\to \cA\times\cA',$$ 
where $\cA$ is an isotrivial abelian scheme with constant fiber $A_0$
and $\omega|_Z\in H^0(\cJ_{\alpha}|_Z,
\Omega^1_{\cJ_{\alpha}|_Z/\nH_\alpha|_Z})$ is the pullback of a
constant 1-form on $\cA$. 
\end{lem}

\begin{rem}
This lemma generalises the construction of \Cref{subsec_Hurwitz} corresponding to the case $\dim A_0 = 1$.
\end{rem}

\begin{proof}[Proof of \Cref{lem_alg_in_isoper}]
Translating the statement in terms of $\QQ$VHS, it is equivalent to prove a splitting
$$\bV_{\alpha,\QQ}|_{Z}\cong \bV'\times\bV'',$$
where $\bV'$ is a constant $\QQ$VHS and such that $\omega|_{Z}$, which
is a global section of $\cV_{\alpha}|_{Z}$, is the pullback of a
constant global section of $\bV'\otimes\cO_Z$.  

The fact the $Z$ is contained in an isoperiodic leaf is equivalent to
the fact that $\omega|_{Z}$ is locally flat for the Gauss-Manin
connection on $\cV_{\alpha}|_{Z}$. Since it is a global section, it is
globally flat. In particular, the largest constant local sub-system,
call it $\bV'$, of $\bV_{\alpha,\QQ}|_{Z}$ is non-zero. By the theorem
of the fixed part, $\bV'$ underlies a polarizable $\QQ$VHS of weight 1
and the orthogonal to $\bV'$ with respect to the polarization gives
the required splitting. 
\end{proof}

\begin{prop}\label{prop_bialg_isoper}
Let $Z\subset \nH_\alpha$ be bi-algebraic subvariety contained in an
isoperiodic leaf. Then $Z$ is affine in the period charts.
\end{prop}

\begin{proof} 
Up to passing to a finite cover of the stratum, we may assume that
there are sections $\sigma_1,\ldots, \sigma_n: \nH_\alpha\to
\cC_\alpha$, which label the zeros, where
$n:=\mathrm{length}(\alpha)$. If $n=1$, then $\nH_\alpha$ is a minimal
stratum and the irreducible components of isoperiodic leaves are
singletons. Henceforth, we assume $n\ge 2$. By
\Cref{lem_alg_in_isoper}, the tautological section $\omega|_Z$ is the
pullback of a constant non-zero 1-form, that we keep calling $\omega$,
on a constant factor $A$ of $\cJ_\alpha|_Z$ of dimension $r\geq 1$. Complete $\omega$ to a
basis $\omega=\omega_1,\ldots,\omega_r$ of
$H^0(A,\Omega^1)$. The dual basis
gives an identification $\Lie(A)\cong \CC^r$, and we denote by
$\pi:\CC^r\to A^\an$ the universal cover, given by the exponential
map. Its multivalued inverse is given by  
\begin{align*}
A^\an &\to \CC^r, \\
x&\mapsto \left( \int_0^x\omega_1,\ldots, \int_0^x\omega_r \right).
\end{align*}
Let $\pr_1:\CC^r\to \CC$ be the projection to the first factor. Consider the following commutative diagram:

\begin{center}
\begin{tikzcd}
\widetilde{Z} \arrow[rr, "\tilde{\sigma}"] \arrow[d] \arrow[rrrr, "f",
bend left]    &  & \CC^r\times\ldots\times\CC^r \arrow[rr, "\pr"]
\arrow[d, "\Pi"] &  & \CC\times\ldots\times \CC \\ 
Z^\an \arrow[rr, "\sigma"] &  & A^\an\times\ldots\times A^\an \,\, .                                                 &  &                        
\end{tikzcd}
\end{center}

where $\pr := (\pr_1,\ldots,\pr_1)$, $\Pi:=(\pi,\ldots,\pi)$,
$\sigma:=(\sigma_2-\sigma_1,\ldots,\sigma_n-\sigma_1)$, and the
products all consist of $n-1$ factors. Let
$V:=\overline{f(\widetilde{Z})}^{\Zar}\subset \CC^{n-1}$ and
$W:=\pr^{-1}(V)$. We want to show that $V$ is an affine subspace. If
$V= \CC^{n-1}$, we are done. Assume $V\subsetneq\CC^{n-1}$, so that
$\codim_{\CC^{r(n-1)}} W\ge 1$. Consider the subvariety
$\sigma(Z)\times W\subset A^{n-1}\times(\CC^r)^{n-1}$, let
$\Delta\subset A^{n-1}\times(\CC^r)^{n-1}$ be the graph of $\Pi$  and
let $U$ be the analytic irreducible component of $(\sigma(Z)\times
W)\cap \Delta$ containing the image of $\widetilde{Z}$. Clearly $\dim
U\ge \dim \sigma(Z)$. The following computation shows that $U$ is an
unlikely intersection (the codimensions are calculated in
$A^{n-1}\times(\CC^r)^{n-1}$ and $N:=r(n-1)$): 

\begin{align*}
&\codim \,\sigma(Z)\times W = 2N-\dim\sigma(Z)-\dim W \ge N-\dim\sigma(Z)+1,\\
&\codim\,\Delta = N,\\
&\codim\,U \le 2N-\dim\sigma(Z).
\end{align*}

According to the Ax-Schanuel theorem (see \Cref{thm_axschanuel}
below), $\pr_{A^{n-1}}(U)$, hence $\sigma(Z)$, is contained in a
translate of a strict abelian subvariety $B$ of $A^{n-1}$. Thus $W$ is
contained in a translate $x+L$ of a strict linear subspace $L\subset
\CC^{r(n-1)}$ and $V\subset \pr(x)+\pr(L)$ (the latter is not
necessarily a strict subset of $\CC^{n-1})$. After translating back to
the origin in $A^{n-1}$, $\CC^{r(n-1)}$ and $\CC^{n-1}$, we obtain a
new commutative diagram 
\begin{center}
\begin{tikzcd}
\widetilde{Z} \arrow[rr, "\tilde{\sigma}"] \arrow[d] \arrow[rrrr, "f",
bend left]    &  & L \arrow[rr, "\pr|_L"] \arrow[d] &  & \pr(L) \\
Z \arrow[rr, "\sigma"] &  & B \,\,
&  &                         
\end{tikzcd}
\end{center}
and we repeat the previous argument. At each step of the process the
dimension of $L$ drops at least by one so the process terminates. In
the final step, the equality $V=\pr(L)$ holds. 
\end{proof}

\begin{theor}[Ax-Schanuel for abelian varieties \cite{Ax72}]\label{thm_axschanuel}
Let $A$ be a complex abelian variety and $\pi:\Lie(A)\to A^\an$ its
universal cover. Let $\Delta\subset A\times \Lie(A)$ be the graph of
$\pi$, $V$ an algebraic subvariety of $A\times\Lie(A)$ and $U$ an
analytic irreducible component of $V\cap \Delta$. If 
$$\codim_{A\times\Lie(A)} U < \codim_{A\times\Lie(A)} V + \codim_{A\times\Lie(A)} \Delta\;\;,$$
then $\pr_A(U)$ is contained in a translate of a strict abelian
subvariety of $A$ (here $\pr_A: A\times \Lie(A)\to A$ is the
projection to the first factor). 
\end{theor}

\begin{proof}[Proof of \Cref{bi-algebraic=linear2}]
Let $Z\subset\nS_\alpha$ be a bi-algebraic subvariety contained in an
isoperiodic leaf and let $Z'$ be the preimage of $Z$ under the
quotient map $\nH_\alpha\to \nS_\alpha$. It is enough to prove that
$Z'$ is linear. 

The model of $Z'$ in the absolute periods is a linear subspace
$L\subset \nV_\alpha^\pure$ of dimension 1, stable under the algebraic
monodromy group $\bfH_{Z'}$ of $\bV_\alpha^\pure|_{Z'}$. Since $\bfH_{Z'}$
is semisimple, it has no non-trivial characters, hence it stabilizes
$L$ pointwise. This implies that the projection of $\omega|_{Z'}$ to
$\cV_\alpha^\pure|_{Z'}$ is stable under the action of $\pi_1(Z')$. By the theorem of
the fixed part, we deduce that $\omega|_{Z'}$ is the pullback of a
(non-constant) 1-form $\omega_A\in H^0(A\times Z', \Omega^1_{A/Z'})$
on a constant factor $A$ of the relative Jacobian
$\cJ_\alpha|_{Z'}$. The assignment $Z'(\CC)\ni z\mapsto
\omega_{A,z}\in H^0(A,\Omega^1)$ is algebraic, and the fibers are
bi-algebraic subvarieties contained in an isoperiodic leaf. By
\Cref{prop_bialg_isoper}, they are affine. Therefore, the algebraic
model $M$ of $Z'$ in the relative periods, is fibered over $L$ by
affine subspaces (i.e. the preimage $M_x$ of $x\in L$ is an affine
subspace of $\nV_\alpha$). But $M$ is stable under $\CC^*$ and
multiplication by $\lambda\in\CC^*$ defines a linear isomorphism
$M_x\to M_{\lambda x}$ for $x\in L$. We conclude that $Z'$ is linear.  
\end{proof}

\subsection{The condition ($\star$) and proof of \Cref{bi-algebraic=linear1}}\label{subsec_cond_star}

\begin{defi} \label{cond_star}
Let $Z\subset\nS_\alpha$ be an irreducible algebraic subvariety. Let
$\mathbb{T}$ be the largest constant local sub-system of $\Gr_1^W
\bV_\alpha|_Z$. We call condition $(\star)$ the following: 
\begin{equation}\label{disp_condition_star}
(\star): \text{There exists $z\in Z(\CC)$ such that $\omega_z\notin \proj\mathbb{T}_z$}. 
\end{equation}
\end{defi}

\begin{proof}[Proof of \Cref{bi-algebraic=linear1}]
We first introduce some notation.
Let $Z\subset \nS_\alpha$ be a bi-algebraic curve satisfying $(\star)$.
We fix a base point $z_0\in Z(\CC)$ corresponding to a pair
$(C_0,[\omega_0])$.
Projection to the pure periods induces a rational map $\pr_\alpha: \proj \nV_\alpha \dashrightarrow \proj \nV_\alpha^\pure$. 
Let $\widetilde{Z}\subset \widetilde{\nS_\alpha^\an}$ be an analytic
irreducible component of $\pi^{-1}(Z^\an)$, where
$\pi:\widetilde{\nS_\alpha^\an}\to \nS_\alpha^\an$ is the universal
cover at $z_0$. 
Let $Y := \overline{D(\widetilde{Z})}^\Zar\subset \proj\nV_\alpha$ be
its algebraic model, where $D$ is the developing map, and
$Y^\pure=\overline{\pr_\alpha(Y)}\subset \proj\nV_\alpha^\pure$ (note that the
image of the developing map does not meet the indeterminacy locus of
$\pr_\alpha$).  
Since $Z$ is bi-algebraic, $Y$ is a closed irreducible curve in 
$\proj\nV_\alpha$.  
Note that $Y^\pure$ is a curve as well: otherwise $Y^\pure$ is a point, thus $Z$ is contained in an isoperiodic leaf, which is not possible since $Z$ satisfies condition $(\star)$.
We need to prove that $Y$ is a linear subvariety. 
We divide the proof in two steps.

\medskip
\noi\textbf{Step 1}: $Y^\pure$ is a linear subvariety. 
\medskip

\noi Let $\bar{v}$ be the image of $z_0$ in $Y^\pure$ and $v$ a lift of $\bar{v}$ to $\nV^\pure_\alpha$.
Let $\bfH$ be the algebraic monodromy group of $\Gr_1^W \bV_{\alpha, \QQ}|_Z$.
Since $\bfH_\CC$ is semisimple, we can uniquely write the $\bfH_\CC$-module $\nV^\pure_\alpha$ as $\nV^\pure_\alpha=T\oplus N$, where $T$ is the largest trivial sub-representation and $N$ its complement. 
Decompose $N$ into its isotypical components 
$$N\cong N_1\oplus\ldots\oplus N_\ell.$$
Observe that the action of $\bfH(\CC)$ on $\proj N_i$ has no fixed points. 
Indeed, let $u\in N_i$, $u\neq 0$ with image $\overline{u}\in\proj N_i$ and assume for the sake of contradiction that $\bfH(\CC)\overline{u}=\overline{u}$. 
Then $\bfH(\CC)$ acts through a character on $\CC u$. 
Since $\bfH_\CC$ is semisimple it has no nontrivial characters, hence $u$ is fixed. 
However the fixed vectors in $\nV^\pure_\alpha$ belong to $T$.

By hypothesis $(\star)$, up to changing base points, we can assume $v\notin T$. 
Write $v=v_T+\sum_{i=1}^\ell v_i$, with $v_T\in T, v_i\in N_i$. Up to
permutation, we can assume that for a certain $\ell'\in\{1,\ldots,
\ell\}$, we have $v_i\neq 0$, for $1\le i\le \ell'$ and $v_i=0$ for
$\ell'<i\le \ell$.  
For $1\le i\le \ell'$, denote $\overline{v_i}$ the image of $v_i$ in $\proj N_i$. 
Then $\dim \bfH(\CC)\overline{v_i}=1$. 
Indeed consider the rational map
\begin{equation*}\label{rat_proj}
p_i:\proj \nV^\pure_\alpha\dashrightarrow \proj N_i. 
\end{equation*}
Since $\overline{v}$ belongs to the domain of definition of this map,
it makes sense to consider the image of $Y^\pure$ under $p_i$ and we have
$\dim p_i(Y^\pure)\le \dim Y^\pure= 1$.   
If $\dim p_i(Y^\pure)=0$ then $\overline{v_i}$ is a fixed point under $\bfH(\CC)$ and we saw that it is not possible. 
We conclude that $\dim p_i(Y^\pure)=1$. 
Since $\bfH(\CC)\overline{v_i}\subset p_i(Y^\pure)$ and $\dim
\bfH(\CC)\overline{v_i}\ge 1$, we deduce $\dim
\bfH(\CC)\overline{v_i}=1$. 
Moreover, the orbit $\bfH(\CC)\overline{v_i}\subset \proj N_i$ is
closed (otherwise, the boundary of the orbit would consist of fixed
points).  
Therefore, the stabilizer of $\overline{v_i}$ is a parabolic subgroup of $\bfH(\CC)$. 
In particular, $\CC v_i$ is stabilized by a Borel subgroup: this can
only happen if $v_i$ is a highest weight vector of the
$\bfH_\CC$-module $N_i$.  

For each $i\in\{1,\ldots,\ell'\}$, let $M_i$ be the
$\bfH_\CC$-submodule of $N_i$ generated by the orbit
$\bfH(\CC)v_i$. It is clearly irreducible and $v_i$ is a
highest-weight vector.  
We now appeal to the classification, due to Satake \cite[Table
p.461]{Sa65} and recalled in \Cref{thm_classif_satake} below, of the
possible irreducible representations that can appear in the
$\bfH_\CC$-module $\nV^\pure_\alpha$.  
We compute the dimension of the orbit of a highest weight vector in each case, following the list of \Cref{thm_classif_satake}. 
The more involved cases are taken care in the representation-theoretic lemma \Cref{lem_reptheory}.
In the end, we will be able to exclude all the cases of Satake's
theorem except $\mathbf{SL}_2$ acting on its standard representation. 
For an algebraic group $\mathbf{G}$ we denote by
$\widetilde{\mathbf{G}}$ its simply connected cover, in the sense of
algebraic groups. 

\begin{enumerate}
\item{Type $A_n$.} Let $\mathbf{SL}_{n+1}$ act on $\bigwedge^r
  \CC^{n+1}$, the $r$-th exterior power of its standard
  representation, with $1\le r\le n$.  
Let $v$ be a highest weight vector of $\bigwedge^r \CC^{n+1}$ and
$\bar{v}$ its image in $\proj \left(\bigwedge^r \CC^{n+1}\right)$.  
Then $\mathbf{SL}_{n+1}(\CC)\bar{v}$ is isomorphic to the Grassmannian
$G(r,n+1)$ parametrizing $r$-planes in $\CC^{n+1}$ and we have $\dim
G(r,n+1)=r(n+1-r)$.  
It has dimension one if and only if $r=1,n=1$, corresponding to $\mathbf{SL}_2$ acting on its standard representation.

\item{Type $B_n$.}  We can assume $n\ge 2$ since $B_1\cong A_1$. Let
  $\widetilde{\mathbf{SO}}_{2n+1}$ act on its spin representation.  
By \Cref{lem_reptheory}(1), the dimension of the orbit of a highest
weight vector is ${n+1 \choose 2}$, which is $\ge 3$ for $n\ge 2$. 

\item{Type $C_n$.} Again we can assume $n\ge 2$ since $C_1\cong A_1$. 
Let $\mathbf{Sp}_{2n}$ act on its standard representation $\CC^{2n}$. 
The orbits are $\{0\}$ and its complement. 
Hence the orbit of a highest weight vector is the whole space and has dimension $\ge 4$.

\item{Type $D_n$.} We can assume $n\ge 3$ because of the isomorphisms
  of Dynkin diagrams $D_1\cong A_1$ and $D_2\cong A_1\times A_1$.  
Let $\widetilde{\mathbf{SO}}_{2n}$ act on one of its half-spin representations. 
By \Cref{lem_reptheory}(2), the dimension of the orbit of a highest
weight vector is ${n \choose 2}$. It never equals $1$ for $n\ge 3$.  
This takes care of the half-spin representations.

The orbits of $\mathbf{SO}_{2n}$ on the standard representation are
the level sets $\{Q(v)=\alpha\}$, where $Q$ is the bilinear form on
$\CC^{2n}$ given by $\begin{pmatrix} 0 & I_n\\ I_n & 0\end{pmatrix}$.  
Since $n\ge 3$, the level sets have complex dimension at least 5.
\end{enumerate}

Write $\widetilde{\bfH}_\CC=\bfH_1\times\cdots\times\bfH_s\times\bfK$
as a direct product where $\bfH_j=\mathbf{SL}_2$ and $\bfK$ is
semisimple without factors isomorphic to $\mathbf{SL}_2$. 
By \Cref{thm_classif_satake} and the previous dimension count, the
action of $\widetilde{\bfH}_\CC$ on $M_i$ must factor through a group
$\bfH_{\sigma(i)}$  where $\sigma:\{1,\ldots, {\ell'}\}\to \{1,\ldots,
s\}$ is a set-theoretic function, and $M_i$ is the standard
representation of $\bfH_{\sigma(i)}$. 
We claim that $\sigma$ maps all indices to the same one: otherwise if
say $\sigma(1)=1, \sigma(2)=2$, then the action of
$\widetilde{\bfH}_\CC$ on $M_1\oplus M_2$ is equivalent to the action
of $\mathbf{SL}_2\times\mathbf{SL}_2$ on $\CC^2\times\CC^2$
coordinatewise, which is transitive.  
This is not possible since the projection of $Y^\pure$ under $\proj
\nV^\pure_\alpha\dashrightarrow \proj (M_1\oplus M_2)$ has dimension one
and is not the whole of $\proj (M_1\oplus M_2)$. 
Finally, since each $M_i$ lives in a different isotypic component, we
deduce that $\ell'=1$ and $\bfH(\CC)v_N=M_1\setminus \{0\}$. 

We now conclude the proof of Step 1. 
Recall that $v=v_T+v_N$, with $v_T\in T$.
If $v_T\neq 0$, then the subset 
$$\{ \lambda(v_T,w): \lambda\in\CC, w\in L\}\subset T\oplus N$$
has dimension 3, and its image in $\proj \nV^\pure_\alpha$, which is
$\bfH(\CC)\overline{v}$, would have dimension $\ge 2$. Therefore
$v_T=0$ and we conclude that $\bfH(\CC)v=M_1\setminus\{0\}$,
i.e. $Y^\pure$ is a projective line.

\medskip
\noi\textbf{Step 2}: $Y$ is a linear subvariety.
\medskip

\noi Let $\bfH^{\rel}$ be the algebraic monodromy group of $\bV_{\alpha, \QQ}|_Z$.
The map $\bfH^\rel\to\bfH$ is surjective.
Let $\cY$ (resp. $\cY^\pure$) the closure of the lifts of $Y$
(resp. $Y^\pure$) to $\nV_\alpha$ (resp. to $\nV_\alpha^\pure$). 
They are closed, two dimensional irreducible subvarieties of their
respective affine spaces. 
By Step 1, $\cY^\pure$ is a linear subspace on
which $\bfH_\CC$ acts through a map $\bfH_\CC\to \mathbf{SL}_{2,\CC}$
and $\cY^\pure$ is isomorphic to the standard representation of
$\mathbf{SL}_{2,\CC}$. 
Note that the complex points of $\bfU:=\ker(\bfH^\rel\to \bfH)^\circ$
act trivially on $\cY$: indeed, the fibers $\cY\to \cY^\pure$ are union of
orbits under $\bfU(\CC)$ and the fibers are zero-dimensional over a
dense open of $\cY^\pure$ (because $\dim\cY=\dim \cY^\pure$). 
Therefore, since $\bfU(\CC)$ is connected, it acts trivially on those fibers, hence everywhere.

Let $W\subset \nV_\alpha$ be the smallest linear subspace containing $\cY$. 
Then $W$ maps surjectively onto $\cY^\pure$. 
Since $\bfU$ acts trivially on $\cY$, it acts trivially on $W$. 
In particular, the action of $\bfH^\rel_\CC$ on $W$ factors through
(the complex points of) $\bfH^\rel/\bfU$, a not necessarily connected
group whose identity component $\bfH'$ is isogenous to $\bfH$. The
group $\bfH'_\CC$ acts on $\cY^\pure$ through $\bfH_\CC$. 
Since $\cY^\pure$ contains a dense orbit under $\bfH(\CC)$ and since $\dim
\cY=\dim\cY^\pure$, it follows that also $\cY$ contains a dense orbit
under $\bfH'(\CC)$.  
This implies that $W$ is an irreducible $\bfH'(\CC)$-module. 
But $W$ surjects onto $\cY^\pure$ hence they are isomorphic as
$\bfH'(\CC)$-modules. In
particular, they have the same dimension. 
This implies $\cY=W$ and that $\bfH^\rel_\CC$ acts on $W$ through the
standard representation of $\mathbf{SL}_{2,\CC}$. 
\end{proof}

\begin{theor}[\cite{Sa65}]\label{thm_classif_satake}
Let $V$ be a $\QQ$-Hodge structure of weight 1, $\mathbf{MT}\subset
\GL(V)$ its Mumford-Tate group and $\bfG$ a connected normal subgroup
of $\mathbf{MT}^{\rm der}$. Let $\widetilde{\bfG}$ be its simply
connected cover and
$\widetilde{\bfG}_\CC=\bfG_1\times\ldots\times\bfG_n$ a decomposition into a product of
almost simple, simply connected $\CC$-groups. Let $W$ be an
irreducible $\bfG_\CC$-submodule of $V_\CC$, and consider it as
$\widetilde{\bfG}_\CC$-module. Write $W\cong W_1\otimes\ldots \otimes W_j$, for some $j\ge 1$, where $W_k$ is an irreducible representation of a factor $\bfG_{i_k}$ of $\widetilde{\bfG}_\CC$, $1\le k\le j$. Then the possible couples $(\bfG_{i_k}, W_k)$, given at the level of the Lie algebras, are as follows
(for $n\ge 1$): 
\begin{itemize}
\item $\mathfrak{sl}_{n+1}(\CC)$, acting on an exterior power $\bigwedge^r \CC^{n+1}$ of the standard representation, $1\le r\le n$.
\item $\mathfrak{so}_{2n+1}(\CC)$, acting on the spin representation, of dimension $2^n$;
\item $\mathfrak{sp}_{2n}(\CC)$, acting on the standard representation, of dimension $2n$;
\item $\mathfrak{so}_{2n}(\CC)$, acting on one of the two half-spin representations, of dimension $2^{n-1}$, or on the standard representation.
\end{itemize}
\end{theor}

\begin{lem}\label{lem_reptheory}
\begin{enumerate}
\item Let $V$ be the spin representation of $\widetilde{\mathbf{SO}}_{2n+1}$, $v\in V$ a highest weight vector and $P$ the parabolic subgroup that stabilizes the line $\CC v$. Then $\dim \widetilde{\mathbf{SO}}_{2n+1}-\dim P={n+1\choose 2}$.
\item Let $V$ be one of the half-spin representations of $\widetilde{\mathbf{SO}}_{2n}$, $v\in V$ a highest weight vector and $P$ the parabolic subgroup that stabilizes the line $\CC v$. Then $\dim \widetilde{\mathbf{SO}}_{2n+1}-\dim P={n \choose 2}$.
\end{enumerate}
\end{lem}
\begin{proof}
Left to the reader.
\end{proof}

%
%
%
%
\subsection{The case of genus 2: proof of \Cref{genus2}}\label{subsec_gen_2}

In genus 2, there are only two strata $\nS_{1,1}$ and $\nS_{2}$, of dimension 4 and 3 respectively.

\begin{prop} \label{6.8}
Let $Z\subset \nS_\alpha$ be a bi-algebraic curve, where
$\alpha=(1,1)$ or $\alpha=(2)$. Then either $Z$ is contained in an
isoperiodic leaf or it satisfies condition $(\star)$. 
\end{prop}
\begin{proof}
Assume that $Z$ does not satisfy $(\star)$. In particular, the fixed
part of $\Gr_1^W \bV_{\alpha,\QQ}|_Z$ is non-zero. There are two
possibilities: $\cJ_\alpha|_Z$ is 
isogenous either to $\cE_1\times \cE_2$, where $\cE_1$ is an isotrivial
elliptic curve over $Z$; or to an isotrivial abelian surface $\cA$ over
$Z$ (and, in these two cases, the tautological section $\omega|_Z$
comes from pullback from the isotrivial factor). In the first case,
$Z$ is contained in the isoperiodic foliation of $\nS_\alpha$ : indeed
an elliptic curve has a unique non-zero 1-form up to
multiplication. (Note that since the isoperiodic foliation on
$\nS_{2}$ has zero-dimensional leaves, we have necessarily
$\alpha=(1,1)$.)  

In the latter case, by Torelli's theorem and the fact that the locus
of abelian varieties isogenous to a fixed abelian variety is
countable, the fibers of 
$\cC_{\alpha}|_Z$ are isomorphic, i.e. $Z$ is contained in a fiber of
the forgetful morphism $\nS_\alpha\to \cM_2$. Since the fibers
$\nS_{2}\to \cM_2$ are zero-dimensional, we must have
$\alpha=(1,1)$. The fibers of $\nS_{1,1}\to \cM_2$ are of dimension
1, so $Z$ actually coincides with a fiber. We prove that such components are not
bi-algebraic. Clearly it is enough to show the corresponding statement
with $\nS_\alpha$ replaced by $\nH_\alpha$ and $\proj \nV_\alpha$ by
$\nV_\alpha$. \Cref{6.8} is thus a consequence of \Cref{prop_vert} below. 
\end{proof}

\begin{prop}\label{prop_vert}
Let $C$ be a smooth projective curve of genus 2. Let $\cU\subset
H^0(C,\Omega^1)$ be the Zariski open subset (complement of a finite
number of lines) consisting of holomorphic 1-forms with two distinct
simple zeros. Consider the multivalued function 
\begin{align*}
T:\cU&\to \CC\\
\omega&\mapsto \int_{\gamma_\omega}\omega
\end{align*}
where $\gamma_\omega\in H_1(C, Z(\omega);\ZZ)$ is a path joining the
distinct zeros, and let $\tilde{T}:\widetilde{\cU}\to \CC$ be a lift
to the universal cover $\pi:\widetilde{\cU}\to\cU$. Then the image of
$(\pi,\widetilde{T}): \widetilde{\cU}\to \cU\times \CC$ is not
``algebraic", i.e. it is not contained in a Zariski closed set of
dimension $=\dim\,\cU=2$. 
\end{prop}

\begin{proof}
Consider the degree two \'etale cover $p:\cU'\to\cU$ on which we can
label the zeros of the 1-form, i.e. there are two morphisms
$\sigma_1,\sigma_2: \cU'\to C$, such that
$\textrm{div}(p(u))=\sigma_1(u)+\sigma_2(u)$, for $u\in\cU'$.  
Consider the $\ZZ$VMHS whose fiber over a point $u$ is
$H_1(C,\{\sigma_1(u),\sigma_2(u)\},\ZZ)$ (it is the dual of the
``restriction" of $\bV_{\alpha,\ZZ}$ to $\cU'$). 
This $\ZZ$VMHS is unipotent. 
Fix a base point $u_0\in \cU'$ and consider the monodromy representation
$$\rho: \pi_1(\cU',u_0)\to\mathrm{GL}(H_1(C,\{\sigma_1(u_0),\sigma_2(u_0)\},\ZZ)).$$
If $\mathrm{im}(\rho)$ is finite then the local system underlying the
$\ZZ$VMHS is constant after a finite \'etale cover. By \cite{HZ} a
unipotent $\ZZ$VMHS is essentially classified by its monodromy
representation, thus our $\ZZ$VMHS is constant. But this contradicts
the remark after (\ref{defi period map}): indeed $\cU$ would map to a point in $\mathfrak{A}_g$ under the period map, contradicting quasi-finiteness.

We deduce that $\mathrm{im}(\rho)$ contains a non-trivial element $A$. If $\gamma$ is
a path connecting $\sigma_1(u_0)$ to $\sigma_2(u_0)$ then $A\gamma =
\gamma + c$, where $c\in H_1(C,\ZZ)$ is non-zero. The set of 1-forms
in $H^0(C,\Omega^1)$ that pair to zero against $c$ is a strict linear
subspace. Let $\omega\in\cU$ be a 1-form that pairs non-trivially
against $c$. Then the fiber of the image of $(\pi,\tilde{T})$ over
$\omega$ contains an infinite coset $x + \ZZ\langle \omega,
c\rangle\subset\CC$ (the number $x\in \CC$ depends on the choice of
the lift $\tilde{T}$). Assume, for the sake of contradiction, that the image of $(\pi,\tilde{T})$ is contained
in a Zariski closed subset $V$ of dimension 2. Since $\dim V = \dim\, \cU$, the (surjective)
projection $V\to \cU$ is generically finite. But the elements
$\omega\in\cU$ with $\langle \omega,c\rangle\neq 0$ are generic in
$\cU$, and by our analysis above, the fiber of $V$ over $\omega$ is infinite. This contradicts generic finiteness. 
\end{proof}

\subsection{An interesting example in genus $3$} \label{genus3}


%


To conclude this paper let us study the bi-algebraic curves in the minimal stratum $\nS_4$, not
satisfying $(\star)$ and not isoperiodic. Let $Z\subset S_\alpha$ be
such a curve. Up to isogeny, the relative Jacobian $\cJ_\alpha|_Z$ factors as
$\cA\times\cA'$ where $\cA$ is an isotrivial abelian scheme over $Z$
corresponding to the fixed part, and $\cA'$ is a complement. By
hypothesis, $\omega|_Z$ comes by pullback from a non-zero relative
one-form on $\cA$. We suppose that $\cA$ is a surface, and call $A$
the abelian surface which is the fiber of $\cA$ above every point.  
We then obtain linearity just by the dimension constraints:

\begin{lem}
$Z$ is linear.
\end{lem}
\begin{proof}[Sketch of proof]
Choose a base point $z_0\in Z(\CC)$ and let $(C_0,[\omega_0])$ be the
corresponding differential. Let $V\subset H^0(C_0,\Omega^1)$ be the
2-dimensional subspace of differentials coming from $A$. Identify $V$
with the corresponding 2-dimensional subspace of $H^1(C_0,\CC)$. Since
$\omega|_Z$ comes by pullback from $A$, its projectivized periods are
constrained on the 1-dimensional line $\proj V\subset \proj
H^1(C_0,\CC)$. 
\end{proof}


Let us now exhibit examples of such curves. The reference here is \cite{HS08}, and we
mostly use their notations. Consider the family of smooth projective
curves parametrized by $\lambda\in\CC - \{ 0,1\}$ given by 
$$C_\lambda : Y^4 = XZ(X-Z)(X-\lambda Z)\subset \proj^2.$$
In the affine chart $\{ Z\neq 0\}$ the curve is given by the equation $y^4=x(x-1)(x-\lambda)$. 
This is the total family of curves that underlies the unique
Shimura-Teichm\"uller curve in genus 3, see \cite{MolST}.

The curves $C_\lambda$ have genus 3 and are not hyperelliptic.
We call $P_0, P_1, P_\lambda, P_\infty$ the four points of $C_\lambda$
lying on the line $\{Y=0\}$, corresponding to $[0,0,1]$, $[1,0,1]$,
$[\lambda,0,1]$ and $[1,0,0]$ respectively. 
On $C_\lambda$, we construct four abelian differentials up to scaling,
with a single zero at $P_i$, $i=0,1,\lambda,\infty$, respectively.  

We first state a trivial but useful fact. Let $C$ be a smooth
projective curve of genus $g$, non-hyperelliptic, and $i:C\to \proj
H^0(C,\Omega^1)^*$ be the canonical embedding. Tautologically,
differential forms up to scaling correspond to hyperplanes in $\proj
H^0(C,\Omega^1)^*$ and the order of vanishing at $x\in C(\CC)$ of a
differential form $\omega$, is the order of tangency at $i(x)$ of
$i(C)$ with the hyperplane $\omega = 0$. 

The embedding in $\proj^2$ coming from the definition of $C_\lambda$
coincides with the canonical embedding (there is only one linear
series $\mathfrak{g}^{g-1}_{2g-2}$ on a given curve of genus $g$). Hence to determine our differentials we
will consider lines in $\proj^2$. 

We collect here the facts that we need on the family $C_\lambda$.

\begin{enumerate}
\item There is an obvious map $f_\lambda$ to the elliptic curve
  $E_{\lambda}=\{ y^2 = x(x-1)(x-\lambda)\}$ given by $(x,y)\mapsto
  (x,y^2)$. 
\item There are involutions $\sigma_\lambda$ and $-\sigma_\lambda$ of
  $C_\lambda$, such that $C_\lambda/\sigma_\lambda \cong C_\lambda/
  {-\sigma_\lambda} \cong E_{-1}$. We call the quotient maps
  $\kappa_{\sigma_\lambda}$ and $\kappa_{-\sigma_\lambda}$
  respectively. [{\em loc. cit.}, Section 1.3 and Proposition 6] 

\item The map  $Jac(C_\lambda)\to E_\lambda\times E_{-1}\times E_{-1}$
  induced by $f_\lambda\times
  \kappa_{\sigma_\lambda}\times\kappa_{-\sigma_\lambda}$ is an
  isogeny. [{\em loc. cit.}, Proposition 7] 
\end{enumerate}
 
We compute the lines corresponding to the unique differential up to
scaling coming from the three factors. This is equivalent to knowing
the ramification points of each of the three maps $f_\lambda,
\kappa_{\sigma_\lambda}, \kappa_{-\sigma_\lambda}$; for $f_\lambda$
they are the four points $P_i$, while for the other two morphisms they
are listed in \cite[Proposition 14]{HS08}.  
Choose a square root $\zeta$ of $1-\lambda$ (here we move to the
degree 2 unramified cover of $\CC-\{0,1\}$ where the square root is
well-defined). 
We find the following:
the differential coming from $E_\lambda$ corresponds to the line $\{ Y = 0 \}$;
the differential coming from $E_{-1}$ by pullback along
$\kappa_{\sigma_\lambda}$ is the line $L_+:=\{ X - (1+ \zeta)Z = 0\}$;
and the differential coming from $E_{-1}$ by pullback along
$\kappa_{-\sigma_\lambda}$ is the line $L_-:=\{ X - (1- \zeta)Z =
0\}$. 
Therefore the differentials up to scaling coming from the factor 
$E_{-1}\times E_{-1}$ correspond to the pencil of lines $\{ \alpha X +
\beta Z = 0 : \alpha,\beta\in\CC \}$, with homogenous coordinates
$[\alpha,\beta]$ (in the chart $\{Z\neq 0\}$ the pencil is the pencil
of vertical lines). 

One easily checks there are exactly four points on $C_\lambda$ such
that the tangent line at the point belongs to the pencil, i.e. the
points $P_0, P_1, P_\lambda, P_\infty$. Moreover at each of these
points the tangent line has order of tangency equal to 4. The
corresponding differentials up to scaling provide our four examples.

\end{document}